%% file: parareal_micmac.tex
\documentclass[a4paper]{svjour3}
\usepackage{geometry}
\usepackage{mathtools}
\mathtoolsset{showonlyrefs}
\usepackage{graphicx}
\usepackage{amsmath,amssymb,amsfonts}
\usepackage{natbib}
\usepackage{tikz}
\usetikzlibrary{calc,shapes,backgrounds,fit,arrows,positioning,decorations.pathmorphing}
\usepackage{pgfplots}
\usepackage{csvsimple}
\usepackage{dcolumn}
\usepackage[utf8]{inputenc}
\usepackage[T1]{fontenc}
\input{macros.tex}

\renewcommand{\vec}[1]{\underline{#1}}

\newcommand{\fineDt}{\mathcal{F}_{\Delta t}}
\newcommand{\coarseDt}{\mathcal{C}_{\Delta t}}
\newcommand{\jump}[1]{\mathcal{J}^{\mathrm{#1}}}
\newcommand{\oper}[3]{\mathcal{#1}}%^{\textrm{#2#3}}}
\newcommand{\coarsesub}{\rho}
\newcommand{\noisysub}{\mathcal{Z}}
\newcommand{\particlesub}{\mathcal{X}}
\newcommand{\finesub}{\mu}

\begin{document}

\title{Parareal computation of stochastic differential equations with time-scale separation: a numerical study}

\author{Fr\'ed\'eric Legoll \and Tony Leli\`evre \and Keith Myerscough \and Giovanni Samaey}

\authorrunning{F. Legoll et al.} % if too long for running head

\institute{F. Legoll \at
              Ecole des Ponts ParisTech and Inria, 6-8 avenue Blaise Pascal, Cit\'e Descartes, 77455 Marne la Vall\'ee, France \\
              \email{frederic.legoll@enpc.fr}           %  \\
           \and
           T. Leli\`evre \at
           Ecole des Ponts ParisTech and Inria, 6-8 avenue Blaise Pascal, Cit\'e Descartes, 77455 Marne la Vall\'ee, France \\
              \email{tony.lelievre@enpc.fr}           %  \\
           \and 
           K. Myerscough \at Department of Computer Science, KU Leuven, Celestijnenlaan 200A, B-3001 Leuven, Belgium \\ \email{keith@myerscough.nl}
           \and 
           G. Samaey \at 
           Department of Computer Science, KU Leuven, Celestijnenlaan 200A, B-3001 Leuven, Belgium\\
           \email{giovanni.samaey@cs.kuleuven.be}           
}

\maketitle

\begin{abstract}
The parareal algorithm is known to allow for a significant reduction in wall clock time for accurate numerical solutions by parallelising across the time dimension. We present and test a micro-macro version of parareal, in which the fine propagator is based on a (high-dimensional, slow-fast) stochastic microscopic model, and the coarse propagator is based on a low-dimensional approximate effective dynamics at slow time scales.
At the microscopic level, we use an ensemble of Monte Carlo particles, whereas the approximate coarse propagator uses the (deterministic) Fokker-Planck equation for the slow degrees of freedom. 
The required coupling between microscopic and macroscopic representations of the system introduces several design options, specifically on how to generate a microscopic probability distribution consistent with a required macroscopic probability distribution and how to perform the coarse-level updating of the macroscopic probability distribution in a meaningful manner.
We numerically study how these design options affect the efficiency of the algorithm in a number of situations. 
The choice of the coarse-level updating operator strongly impacts the result, with a superior performance if addition and subtraction of the quantile function (inverse cumulative distribution) is used. 
How microscopic states are generated has a less pronounced impact, provided a suitable prior microscopic state is used. 
\end{abstract}

\maketitle

\section{Introduction}

In many applications, a system is modelled using a high-dimensional system of stochastic differential equations (SDEs) that capture phenomena occurring at multiple time scales. 
Usually the quantities of interest are \emph{macroscopic} observables of the system, lower-dimensional than the full (\emph{microscopic}) state, and evolving on a slower time scale than the fine-scale, microscopic dynamics. In such a setting, the computational cost of simulating the microscopic system over macroscopic time intervals may be prohibitive, due to both time step restrictions and the number of degrees of freedom. 
One then often resorts to simulations using low-dimensional (coarse-grained, effective) models, in which the fast degrees of freedom are eliminated. 
Typically, this is done by taking a limit in which the time-scale separation is infinite, which corresponds to the assumption that the (eliminated) fast degrees of freedom equilibrate infinitely quickly with respect to the (retained) slow ones. 
\citet{PaSt08}, for example, present a review of proposed methods to obtain such macroscopic models, either analytically or numerically. See also \citep{LeLe10} and subsequent works.
Several classes of numerical approaches have been proposed to address this type of problems: we refer, for instance, to the works on equation-free \citep{KeGeHyKeRuTh03, KeSa09} or heterogeneous multiscale methods \citep{EEn03, EEnLiReVa07}, and references therein.

Approximate macroscopic models perform well when the time-scale separation is large and the approximation of an infinite time-scale separation is acceptable. 
For moderate time-scale separation, however, the macroscopic model is usually insufficiently accurate. Simultaneously, directly simulating the microscopic dynamics may still be prohibitively expensive. 
As a consequence, there is a need for computational multiscale methods which converge, for a limited CPU cost, to the full microscopic dynamics, even in cases when the time-scale separation is fixed and finite. 

In this article, we present and test a micro-macro version of the \emph{parareal} algorithm of \citet{LiMaTu01} for high-dimensional scale-separated SDEs. Parareal was introduced to solve time-dependent problems using computations in parallel, exploiting the presence of multiple processors to reduce the real (wall-clock) time needed to obtain a solution on a long time interval. 
It is based on a decomposition of the time interval into subintervals, and makes use of a predictor-corrector strategy: a coarse (cheap, but inaccurate) propagator is used to obtain a prediction that is iteratively corrected using a fine (expensive, but accurate) propagator. 
The calculation of the corrections is performed in parallel across the time domain, concurrently on the different available processors. 
If a sufficiently accurate result is obtained in only a few iterations, compared to the number of subintervals, a large reduction of the required wall-clock time is achieved with respect to a full microscopic simulation. 
The parareal algorithm has previously found application in many different fields, such as nuclear fusion \citep{SaNeSa10}, power grid systems \citep{GuDiSrSiSt15}, the automotive industry \citep{LoHeLo14}, molecular dynamics \citep{PaMa14}, $N$-body problems \citep{SpRuKrEmMiWiGi12}, and medical applications \citep{RaKa14,KrNaRuSpWiKr15}, to mention but a few.

A micro-macro version of parareal was proposed by \citep{LeLeSa13} for singularly perturbed systems of ordinary differential equations (ODEs). It uses a coarse propagator based on an approximate macroscopic model with fewer degrees of freedom. Here, we construct a micro-macro parareal method in which the fine propagator is a Monte Carlo discretization of a high-dimensional slow-fast SDE. The coarse propagator uses as the approximate macroscopic model a limiting (low-dimensional) Fokker-Planck equation for the slow degrees of freedom that arises in the limit of infinite time-scale separation. Two main additional difficulties arise in this setting, compared to the setting of \citep{LeLeSa13}. 
First, we need to define an operator that creates a microscopic probability distribution consistent with an imposed marginal distribution of the slow degrees of freedom, and draw samples from this distribution. We seek the desired microscopic state as a minimal perturbation of a previously available (prior) microscopic state; a procedure that we call \emph{matching}. Matching methods have also been developed for scale-separated SDEs, with the macroscopic state being a number of moments of the microscopic probability distribution, see \citep{DeSaZi15} and \citep{LeSaZi17}.
Second, we need to propose a meaningful parareal correction strategy when the system state is a discretized probability distribution. 
Classical parareal \citep{LiMaTu01} applies these updates in an additive way. 
This procedure is not appropriate for the Fokker-Planck setting, since the modelled quantities are probability distributions, which need to remain positive and of unit mass. 
Both for the matching and the parareal corrections, there are several options and it is not a priori clear how these options influence the performance of the micro-macro parareal method. In this article, we therefore present a generic micro-macro parareal algorithm for scale-separated SDEs, with a detailed numerical comparison of the performance of the algorithm, for various choices for both procedures. 

While the micro-macro parareal method that we present in this work converges to the full microscopic dynamics at a reduced wall-clock time, the total CPU time that is consumed (summed over all processors) is larger than that for one direct microscopic simulation, due to the iterative character of the method. 
Nevertheless, there are several scenarios in which parareal may be computationally beneficial. 
A first situation is when the number of available processors is so large that they cannot be put to work sufficiently by other means of parallelisation, such as parallelising the Monte Carlo realisations. This happens, for instance, when we are interested in simulations over extremely long time horizons. 
Now that the developments in computer hardware focus on increasing the number of cores on a CPU, and no longer on increasing the clock speed, we expect such situations to occur more frequently. 
A second possible setting is that of PDE-constrained optimisation. 
In that setting, one needs to solve a forward and backward problem simultaneously, which needs to be done iteratively anyway, even when the PDEs are solved using time-marching. 
Then, parallelisation in time does not introduce an extra iteration (see e.g. the work of \cite{guenther2018pintoptim}).
Finally, in real-time model predictive control problems, the dominant constraint on simulations is on the wall-clock time, since simulations must be finished for their output to be used in control decisions. 
There is then a clear interest in reducing the wall-clock time for a given simulation, even at the cost of increasing the total required CPU time. 
A real-life application that combines the need for real-time control with simulations over long time horizons is the chemical process industry, in which control actions on the chemical reactor may need to be taken several hours in advance.  Since chemical reaction kinetics is a prototypical application of atomistic modelling techniques, the model problem in the present manuscript is also drawn from this particular field of application. 

The remainder of this article is organised as follows. 
Section~\ref{sec:micmac} discusses the model problem that we use throughout the article. 
We explain the generic structure of the micro-macro parareal method in Section~\ref{sec:parareal}. We then turn to the two main building blocks that define the particular method: the computation of the parareal corrections (Section~\ref{sec:jumps}) and the matching operator (Section~\ref{sec:matching}).
Section~\ref{sec:results} then contains extensive numerical experiments, comparing the performance of the different options. 
We conclude in Section~\ref{sec:conclusion}.

\section{Micro-macro modeling}
\label{sec:micmac}

This section presents the slow-fast stochastic differential equation (SDE) that serves as the microscopic model (Section~\ref{sec:mic}) and the approximate macroscopic Fokker-Planck equation derived from it under the assumption that the fast degrees of freedom have equilibrated (Section~\ref{sec:mac}). 

\subsection{Microscopic model}
\label{sec:mic}

\subsubsection{Slow-fast SDEs} 

As a microscopic model, we consider a slow-fast system of coupled It\=o SDEs for $X_t \in \IR$ and $Y_t \in \IR^d$, which we write in a general form (see, e.g., \citet{GiKuSt04}) as
\begin{align}
\dd X_t &= \frac{-1}{\tau_x} \partial_x F(X_t, Y_t) \dd t + \sqrt{\frac{2 D(X_t, Y_t)}{\tau_x}} \dd U_t \label{eq:SDEx}\\
\dd Y_t &= \frac{-1}{\epsilon} \partial_y G(X_t, Y_t) \dd t + \sqrt{\frac{2 E(X_t, Y_t)}{\epsilon}} \dd V_t \label{eq:SDEy},
\end{align}
where $U_t$ and $V_t$ denote two independent Wiener processes. 
In the numerical examples of this article, the fast variable $Y$ is a scalar, but all algorithms apply straightforwardly when the fast variables are higher-dimensional.
Initial conditions are specified as the joint probability distribution for $X_0$ and $Y_0$. 
These will be given later for the cases discussed in this article, but will always be readily sampled distributions, such as independent lognormal distributions.

We denote by $\mu(t,x,y)$ the particle distribution, i.e., we have
\[
P\left[ (X_t, Y_t) \in A \right] = \int_A \mu(t, x, y) \dd x \dd y.
\] 
The evolution of this distribution is governed by the Fokker-Planck equation
%\begin{align}
%%\mu_t(t,x,y) 
%%&= \mu(t,x,y) F_{xx}(x,y) + \mu_x(t,x,y) F_x(x,y) + \frac{1}{\epsilon} \mu(t,x,y) G_{yy}(x,y) +\\
%%&+ \frac{1}{\epsilon} \mu_y(t,x,y) G_y(x,y) + \frac{1}{\beta} \mu_{xx}(t,x,y) + \frac{1}{\beta\epsilon} \mu_{yy}(t,x,y),\label{eq:FPxy}
%\partial_t \mu(t,x,y) 
%&= \partial_x \left( \mu(t,x,y) \partial_x F(x,y) \right) + 
%\frac{1}{\epsilon} \partial_y \left( \mu(t,x,y) \partial_y G(x,y) \right) \\
%&+ 
%\frac{1}{\beta} \partial_x^2 \mu(t,x,y) + \frac{1}{\beta\epsilon} \partial_y^2 \mu(t,x,y).\label{eq:FPxy}
%\end{align}
\begin{align}
\partial_t \mu = \frac{1}{\tau_x} \partial_x \left( \mu \partial_x F + \partial_x \left( \mu D \right)\right) + 
\frac{1}{\epsilon} \partial_y \left( \mu \partial_y G + \partial_y \left( \mu E \right)\right). \label{eq:FPxy}
\end{align}

We choose the parameters $\tau_x\gg\epsilon>0$ to ensure a time-scale separation between the (slow) evolution of $X_t$ and the (fast) evolution of $Y_t$. 

\subsubsection{Example: a chemically reacting system}
\label{sec:bruna_mic}

Motivated by the above discussion, we consider as the microscopic model a slow-fast SDE that models chemical reactions, taken from \citet{BrChSm14}, which is given by~\eqref{eq:SDEx}--\eqref{eq:SDEy} for a particular choice of potentials and diffusions. The microscopic potentials have the form
\begin{subequations}
\begin{align}
F(x,y) &= -k_4xy - \half k_5x^2, \label{eq:F_Bruna}\\ 
G(x,y) &= -k_0y - \half k_1y^2 + k_2y^2\left(\frac{y}{3}-\frac{1}{2}\right) - k_3y^2\left(\half y-1\right)^2 \label{eq:G_Bruna}.
\end{align} 
\end{subequations}
The diffusion terms are given by
\begin{subequations}
\begin{align}
D(x,y) &= \half\left(k_4 y + k_5 x\right), %\label{eq:D_Bruna}
\nonumber \\ 
E(y) &= \half\Big(k_0 + k_1 y + k_2 y(y-1) + k_3 y(y-1)(y-2)\Big) \label{eq:E_Bruna}.
\end{align} 
\end{subequations}
This is a challenging system due to the presence of bistability of the fast variable for certain parameter choices. 
The bistability introduces a third (after the fast and slow time scales $\epsilon$ and $\tau_x$) time scale $\tau_s$, that indicates the switching time between the two wells of the fast variable. 
The behaviour of the system is different depending on the relative time scales.
To be precise, \citet{BrChSm14} identify three different regimes:
\begin{itemize}
\item \emph{Regime 1 (well equilibrated)} is characterized by a fast switching time $\tau_s$ of the fast variable compared to the time scale of the slow variable $\tau_x$, i.e.~$\epsilon \ll \tau_s \ll \tau_x$. 
\item \emph{Regime 2 (metastable fast dynamics)} corresponds to a balanced switching time and slow variable evolution, i.e.~$\epsilon \ll \tau_s \sim \tau_x$.
\item \emph{Regime 3 (no fast dynamics)} implies the slow variable evolves at a similar speed as the fast variable, but the fast variable has a slow switching time, i.e.~$\epsilon \sim \tau_x \ll \tau_s$. 
\end{itemize}
In this article, we are only interested in ``Regime 1'' and ``Regime 2''. 

\subsubsection{Weighted Monte Carlo simulation}

We use weighted Monte Carlo samples evolved by the Euler-Maruyama scheme as the fine propagator in our multiscale simulations. 
This choice is motivated by the goal of (eventually) simulating systems with a higher dimensional fast component, where a grid-based discretization of~\eqref{eq:FPxy} would suffer from the \emph{curse of dimensionality} \citep{Bellman57}. 

The SDEs~\eqref{eq:SDEx}--\eqref{eq:SDEy} can be represented by a finite number of samples $(X_t,Y_t)$, possibly with attached weights $W_t$. 
Such samples are gathered in an ensemble
\begin{align*}
\mathcal{X}_t = \ensemble{X_t^p, Y_t^p, W_t^p}{p=0,1,\dots,P-1}.
\end{align*}
The density sampled by the microscopic ensemble at a given time $t$ is given by $\mu(t,x,y)$ such that
\begin{align}
\forall A \subset \IR^2, \quad \int_A \mu(t,x,y) \dd x \dd y = \expbra{\sum_{p \text{ s.t. } (X_t^p, Y_t^p) \in A} W_t^p}.
\label{eq:MCsample}
\end{align}
Note that we choose to normalise the weights such that $\sum_{p=0}^{P-1} W^p_t = 1$ at any time $t$, thus there is no need for a normalising quotient in~\eqref{eq:MCsample}. The expectation for a certain observable $M(x,y)$ of the stochastic process at time $t$ is given by the Monte-Carlo estimator
\begin{align}
\widehat{M}_t = \sum_{p=0}^{P-1} M(X_t^p, Y_t^p) \, W_t^p.
\label{eq:MCestimator}
\end{align}
The estimator~\eqref{eq:MCestimator} is unbiased, that is
\begin{align*}
\expbra{\widehat{M}_t} = \int M(x,y) \, \mu(t,x,y) \dd x \dd y. 
\end{align*}

We simulate the evolution of the ensemble using the Euler-Maruyama scheme, using a time step $\delta t$ and approximating the system state $(X_{m\delta t}^p, Y_{m\delta t}^p)$ by the tuple $(X_m^p, Y_m^p)$.
Beware that there is a slight abuse of notation, with $X$ and $Y$ denoting both the exact solution to the SDEs~\eqref{eq:SDEx}--\eqref{eq:SDEy} as well as the approximation; the distinction being marked by the subscript index. 
The weights are not used in determining the evolution of the ensemble, but only appear when estimating observables via~\eqref{eq:MCestimator}. 
The Euler-Maruyama scheme reads as
\begin{subequations}
\begin{align}
X_{m+1}^p = X_m^p - 
\frac{1}{\tau_x} \partial_x F(X_m^p, Y_m^p) \, \delta t + 
\sqrt{\frac{2 D(X_m^p, Y_m^p)}{\tau_x}} \, \delta U_m, \label{eq:EMx} \\
Y_{m+1}^p = Y_m^p - 
\frac{1}{\epsilon} \partial_y G(X_m^p, Y_m^p) \, \delta t + 
\sqrt{\frac{2 E(X_m^p, Y_m^p)}{\epsilon}} \, \delta V_m, \label{eq:EMy}
\end{align}
\end{subequations}
where $\delta U_m$ and $\delta V_m$ are (pseudo-)random numbers drawn independently from the normal distribution $\mathcal{N}(0,\delta t)$ with mean zero and variance $\delta t$. 
Again, beware that there is a slight abuse of notation, with $U$ and $V$ representing discrete analogues to the Brownian motions represented by the same symbols in~\eqref{eq:SDEx}--\eqref{eq:SDEy}. Also here, the distinction should be clear from the used subscripts. 

\begin{remark}
At this point, we do not specify any evolution law for the weights $W$. They can, for now, be considered to be determined when discretising the initial condition via a weighted particle ensemble, and kept constant during evolution. In Section~\ref{sec:matching}, we will use the weights to implement the parareal parts of the algorithm.
\end{remark}

The microscopic simulation is costly for two reasons. First, the addition of the fast degrees of freedom increases the dimensionality of the evolving probability distribution, and hence, the number of unknowns that need to be tracked as a function of time. Second, the presence of fast time scales severely limits the maximal allowed time step $\delta t$. In particular, for the slow-fast system~\eqref{eq:SDEx}--\eqref{eq:SDEy}, we are forced to choose $\delta t=\mathcal{O}(\epsilon)$, which may be prohibitively small.

In the remainder of this article, we assume that the result of the approximation~\eqref{eq:EMx}--\eqref{eq:EMy} is truthful, i.e., we neglect time discretization errors and only consider a statistical error at the microscopic level due to the finite size of the ensemble $\mathcal{X}$.  We are \emph{not} concerned with improving the fidelity of this microscopic propagator. 

%For completeness we do propose a discrete represenation of the two-dimensional probability density $\mu(x,y,t)$ on a regular grid with $x_0 < x_1 < \dots < x_Q$ and $y_0 < y_1 < \dots < y_Q$. 
%We use the integral quantities 
%\begin{align}
%\vec\mu(t) = \ensemble{\mu_{ij}(t)}{i,j=0,1,\dots,Q-1}, \mbox{ where each } \mu_{ij}(t) = \int_{y_i}^{y_i+1} \int_{x_i}^{x_i+1} \mu(x,y,t) \dd x \dd y,
%\end{align}
%which we will define to be elements of the phase space $D_\particlesub = \IR^{Q^2}$. 

\subsection{Macroscopic model}
\label{sec:mac}

\subsubsection{Macroscopic Fokker-Planck equation}

In the limit of infinite time scale separation, that is when $\epsilon \searrow 0$, an approximate macroscopic model can be obtained in order to describe the evolution of only the slow variable $X_t$: 
\begin{align}
\dd Z_t &= \frac{-1}{\tau_x}\partial_x \overline{F}(Z_t) + \sqrt{\frac{2 \overline{D}(Z_t)}{\tau_x}} \dd W_t^x, \label{eq:SDEz}
\end{align}
where $W_t^x$ is a Wiener process and where the effective potential $\overline{F}(x)$ is such that its derivative is equal to the expectation of $\partial_x F(x,y)$ with respect to invariant distribution $\phi(y;x)$ of the fast variable $y$, given the slow variable $x$, i.e., the invariant measure of equation~\eqref{eq:SDEy} for a fixed and given value of $x$, see \citep{PaSt08}:
\begin{align}
%\partial_x \overline{F}(x) = \frac{ \int \partial_x F(x,y) \phi(y;x) \dd y }{ \int \phi(y;x) \dd y }.
\partial_x \overline{F}(x) = \int \partial_x F(x,y) \, \phi(y;x) \dd y.
\label{eq:Fbar}
\end{align}
The effective diffusion is found in a similar manner to be
\begin{align}
%\overline{D}(x) = \frac{ \int D(x,y) \phi(y;x) \dd y }{ \int \phi(y;x) \dd y }.
\overline{D}(x) =  \int D(x,y) \, \phi(y;x) \dd y.
\label{eq:Dbar}
\end{align}
Associated with~\eqref{eq:SDEz} is the Fokker-Planck equation %\todo{Corrected equation below - check code for what is done.}
%\begin{align}
%\partial_t \rho = \frac{1}{\tau_x} \partial_x \left( - \rho \partial_x \overline{F}_x  + 2  \partial_x (\overline{D}(x) %\rho)\right).\label{eq:FPz}
%\end{align}
\begin{align}
\partial_t \rho = \frac{1}{\tau_x} \partial_x \big( \rho \partial_x \overline{F} + \partial_x (\overline{D} \rho)\big).\label{eq:FPz}
\end{align}

%\begin{align}
%\partial_t \rho(t,x) = \partial_x \left( \rho(t,x) \partial_x \overline{F}_x(x) \right) + \frac{1}{\beta} \partial_x^2 \rho(t,x).\label{eq:FPz}
%\end{align}

\subsubsection{Finite volume scheme}\label{sec:mac_fv}

The parareal algorithm is only a sensible tool for accelerating simulations if an approximate simulation using the coarse propagator is much cheaper than a full simulation using the microscopic model. 
This is possible here because the macroscopic system, on which the coarse propagator is based, is low-dimensional. Because of this low-dimensionality, we are able to use a deterministic, grid-based discretisation of the Fokker-Planck equation~\eqref{eq:FPz}, rather than a Monte-Carlo simulation as was used for the microscopic system. %, as this is computationally cheaper for the same accuracy. 
In the discretization below, we additionally take advantage of the one-dimensionality of the macroscopic Fokker-Planck equation.

The density $\rho$ is discretised in space on a grid $x_0 < x_1 < \dots < x_I$ to yield the vector
\begin{align}
\vec{\rho}(t) = \ensemble{\rho_i(t)}{i=0,1,\dots,I-1} \quad \mbox{with} \quad \rho_i(t) = \int_{x_i}^{x_{i+1}} \rho(t,x) \dd x,
\label{eq:integral}
\end{align}
which is an element of the phase space $D_{\coarsesub} = \IR^I$.
In the remainder of this article, we will often view the discretized probability density $\vec{\rho}(t)$ as a piecewise constant function
\begin{align}
\rho_\mathrm{pc}(x) = \sum_{i=0}^{I-1} \frac{\rho_i}{x_{i+1} - x_i} \ \bs{1}_{(x_i,x_{i+1}]}(x).
\label{eq:piecewise_constant_pdf}
\end{align}

To determine the evolution of $\vec{\rho}(t)$, we first rewrite~\eqref{eq:FPz} as
\begin{align}
\partial_t \rho = \frac{1}{\tau_x} \partial_x \left( \overline{D} \, \rho^\infty \, \partial_x \left( \frac{\rho}{\rho^\infty} \right) \right),\label{eq:FPz_rewritten}
\end{align}
where 
\[
\rho^\infty(x) = \frac{Z^\infty}{\overline{D}(x)} \exp\left(- \int_0^x \frac{\partial_x \overline{F}(s)}{\overline{D}(s)} \dd s \right), 
\]
 and $Z^\infty$ is a normalising coefficient that need not be computed. %\footnote{This coefficient need not be computed, as it is divided out because $\rho^\infty$ appears in~\eqref{eq:FPz_rewritten} once in inverted form.}. 
Integrating~\eqref{eq:FPz_rewritten} over $[x_i, x_{i+1}]$ yields
\begin{align*}
\partial_t \rho_i(t) = \frac{1}{\tau_x} \big( \phi_{i+1}(t) - \phi_i(t) \big),
\end{align*}
where the fluxes are given by $\displaystyle \phi_i(t) = \overline{D}(x_i) \, \rho^\infty(x_i) \ \partial_x \left( \frac{\rho(t, \cdot)}{\rho^\infty} \right) (t, x_i)$.
The remaining derivative is approximated using a central difference scheme:
\begin{align*}
\partial_x \left( \frac{\rho(t, \cdot)}{\rho^\infty} \right) (t, x_i) \approx 
\frac{\rho_{i+1/2}(t)}{h_i \, h_{i+1/2} \, \rho^\infty_{i+1/2}} - \frac{\rho_{i-1/2}(t)}{h_i \, h_{i-1/2} \, \rho^\infty_{i-1/2}},
\end{align*}
where $h_{i+1/2} = x_{i+1} - x_i$, $\displaystyle h_i = \frac{x_{i+1} - x_{i-1}}{2}$ and $\displaystyle \rho^\infty_{i+1/2} = \rho^\infty\left(\frac{x_{i+1} + x_i}{2}\right)$.
Successive substitution results in a linear system that we solve using matrix exponentiation to obtain an exact time integrator for the spatially discretized system $\vec{\rho}(t) \in \IR^I$.
%As such the coarse propagator contains only a spatial discretization error besides the modeling error. 

\subsubsection{The chemical reactions example}

For the specific example of Section~\ref{sec:bruna_mic}, the macroscopic model described above is fairly straightforward. 
Since there are only first order terms in $y$ in the slow potential $F(x,y)$ defined in~\eqref{eq:F_Bruna}, the effective potential $\overline{F}$ is computed using only the mean $\overline{y}$ of the fast variable:
\begin{align*}
\overline{F}(x) &= -k_4 \, x \, \overline{y} - \half k_5 x^2. %\label{eq:Fbar_Bruna}
\end{align*}
Additionally, because the fast potential~\eqref{eq:G_Bruna} and diffusion~\eqref{eq:E_Bruna} do not depend on the slow variable $x$, the mean $\overline{y}$ can be computed once and for all via numerical integration (using the trapezoidal rule) of the equilibrium distribution
\begin{align*}
\phi(y;x) = \phi(y) \sim \frac{1}{E(y)} \exp\left[ -\int_0^y \frac{G'(s)}{E(s)} \dd s \right].
\end{align*} 
Similarly the effective diffusion $\overline{D}(x)$ is given by
\begin{align*} 
\overline{D}(x) &= \half\left(k_4 \overline{y} + k_5 x\right). %\label{eq:Dbar_Bruna}
\end{align*} 

This macroscopic model is an accurate approximation of the dynamics when the coefficients~\eqref{eq:Fbar} and~\eqref{eq:Dbar} are reasonable approximations of the macroscopic drift and diffusion. This is the case when the fast time-scale is such that the full microscopic conditional distribution~$\phi(y;x)$ is reached very quickly, i.e., in ``Regime 1'' as defined in Section~\ref{sec:bruna_mic}. In ``Regime 2'', such an approximation is not justified, because the metastability implies that one only samples around one of the two metastable states on the time scale of the evolution of the slow variable. It is in this regime that the micro-macro parareal method will turn out to deliver a computational advantage.
 
%\begin{figure}
%\centering
%\input{couplings.tikz}
%\caption{Diagram illustrating the different system descriptions and the coupling operators between them.}
%\label{fig:couplings}
%\end{figure}
%

\section{The parareal algorithm}
\label{sec:parareal}

In this section, we introduce the generic structure of the micro-macro parareal algorithm. We recall the traditional parareal algorithm for a generic initial value problem in Section~\ref{sec:parareal_traditional}. We then introduce the micro-macro parareal algorithm in Section~\ref{sec:parareal_micmac}. We discuss some crucial properties of the algorithm and its components in Section~\ref{sec:couplings:requirements}. Specific choices for each of the method's building blocks are deferred to Section~\ref{sec:jumps} and Section~\ref{sec:matching}.

\subsection{The traditional parareal algorithm}
\label{sec:parareal_traditional}

Consider the initial value problem 
\begin{align*}
\fd{u}{t} = f(u), \quad u(0)=u_0, \quad u(t) \in \IR^d, \quad t \in [0,T],
%\label{eq:ivp}
\end{align*}
where $u(t)$ denotes the solution at time $t$, starting from an initial condition $u_0$. 
At each time $t$, the state $u(t)$ is an element of the $d$-dimensional space $\IR^d$ and has an exact flow given by $u(t)$. 
We are interested in finding out the evolution of the system over some time interval $[0,T]$. 
The parareal algorithm divides this time interval $[0,T]$ into $N$ subintervals of length $\Delta t = T/N$ and then iteratively generates a sequence of $N$-tuples $\vec{u}_k \equiv \left\lbrace u_k^n \right\rbrace_{1 \leq n \leq N}$, in which each $u_k^n$ is an approximation to $u(n\Delta t)$ at iteration $k \geq 0$. 
The algorithm requires two numerical propagators that approximate the exact flow over a time interval of length $\Delta t$: (i) a fine propagator $\fineDt$ that is accurate, but requires large computing time; and (ii) a coarse propagator $\coarseDt$ that is less accurate, but is much faster to compute. 
Note that the length of the ``reporting interval'' $\Delta t$ is not necessarily equal to the time step used by either propagator internally. 
In fact, in some implementations, the only difference between $\fineDt$ and $\coarseDt$ is the internal time step. 
The initial approximation, for $k=0$, is given using only the coarse propagator
\begin{align*}
u_{k=0}^{n+1} = \coarseDt (u_{k=0}^n) \quad \mbox{with} \quad u_{k=0}^0 = u_0.
%\label{eq:initial}
\end{align*}
Subsequently, for $k > 0$, the approximation is iteratively improved. To this end, we compute, at each iteration $k$, the corrections
\begin{align*}
 \mathcal{T}_{k+1}^{n+1} := \fineDt (u_k^n) - \coarseDt (u_k^n),
\end{align*}
which can be done efficiently, in terms of wall clock time, \emph{in parallel} for $0 \leq n \leq N-1$. These corrections are added, sequentially, to the result of the coarse propagator,
\begin{align}
u_{k+1}^{n+1} = \coarseDt (u_{k+1}^n) + \mathcal{T}_{k+1}^{n+1},
\label{eq:sequential}
\end{align}
with initial condition $u_{k+1}^{0} = u_0$.
It has been proven that this algorithm converges to the reference solution, given by using the fine propagator sequentially, as $k$ tends to infinity \citep{LiMaTu01}. 

Assuming that the coarse propagator is computationally inexpensive, the computational cost is dominated by the fine propagator. 
The computational time for the sequential propagation of the fine propagator -- the reference solution -- is proportional to $N \Delta t$. 
For the parareal algorithm, the wall-clock time is proportional to $K\Delta t$, depending on the number of required iterations $K$. 
Neglecting the computational overhead of parallelising and the time required for the sequential coarse propagator, the parareal algorithm could achieve a speed-up of $N/K$. 
One can readily show that, for all $n \leq k$, the parareal method is equal to a fine propagation, i.e., $u_k^n = \fineDt^n u_0$. This implies that the required number of iterations is \emph{at most} $N$ to get full convergence over the $N$ time steps. Ignoring parallelisation overhead, the required wall clock time of parareal is never larger than that of naive simulations using the fine propagator. In practice, however, the required number of iterations $K$ should be much smaller than $N$ for the method to be useful, and this is often the case (see for instance \citep{DaLeLeMa13}, resp. \citep{LeLeSa13}, for convergence studies for a Hamiltonian ODE, resp. a two-scale deterministic system).
% We will demonstrate that this property is maintained in the micro-macro setting.
For long simulations that require only a few parareal iterations, the gain in wall-clock time may be substantial. 

\subsection{The micro-macro parareal algorithm}
\label{sec:parareal_micmac}

We now present a micro-macro parareal algorithm that is able to deal with the probabilistic multiscale setting described in Section~\ref{sec:micmac}. Two modifications are required with respect to the standard parareal method of Section~\ref{sec:parareal_traditional}.  First, we intend the coarse (resp.~fine) propagator to be based on a macroscopic (resp.~microscopic) model of the system. Since this implies a different number of degrees of freedom at the two levels, an additional \emph{matching} operator is required to couple the two levels of description. Second, we also need to ensure that the coarse system state, a discretised probability density function, remains positive and of unit mass. This requires changing the additive correction formula~\eqref{eq:sequential}. 

In Section~\ref{sec:prop}, we first describe the coarse and fine propagators. We continue with the operators that connect the microscopic and macroscopic levels of description in Section~\ref{sec:match_restrict}.
Subsequently, we describe the micro-macro parareal algorithm for a generic iteration ($k\ge 1$) in Section~\ref{sec:micmac_parareal_generic}. Afterwards, in Section~\ref{sec:micmac_parareal_zero}, we discuss the zeroth iteration in some more details. 

\subsubsection{Propagators}\label{sec:prop}

For the microscopic model, we consider the high-dimensional slow-fast Fokker-Planck equation~\eqref{eq:FPxy}, which we discretise using a weighted Monte Carlo method.  Thus, as the microscopic state, we introduce an ensemble of $P$ weighted particles, 
\[
\mathcal{X}=\ensemble{X^p, Y^p, W^p}{p=0,1,\dots,P-1}.
\]
Each particle evolves according to the SDE~\eqref{eq:SDEx}--\eqref{eq:SDEy}. For the time discretisation, we use the Euler-Maruyama method~\eqref{eq:EMx}--\eqref{eq:EMy}.  We denote the corresponding fine propagator as 
\begin{align*}
\mathcal{X}^{n+1} = \fineDt (\mathcal{X}^n).
%\label{eq:fine_propagator}
\end{align*}

For the macroscopic model, we use the low-dimensional Fokker-Planck equation~\eqref{eq:FPz}, discretised with the finite volume method described in Section~\ref{sec:mac_fv}.  We denote the coarse propagator as 
\begin{align*}
%\label{eq:coarse_propagator}
\vec\rho^{n+1} = \coarseDt (\vec\rho^n) .
\end{align*}

\subsubsection{Restriction and matching}\label{sec:match_restrict}

\paragraph{Restriction.} Transferring information from the microscopic to the macroscopic level is relatively straightforward. To this end, we introduce a \emph{restriction} operator $\oper{R}{}{}$. The restriction operator maps a microscopic state, a weighted ensemble $\mathcal{X}$, to the corresponding macroscopic state, a discretised probability distribution $\vec\rho$ of only the slow degree of freedom.
We view restriction as a composition of two operators $\mathcal{R} = \mathcal{R}_2 \circ \mathcal{R}_1$.
The first operator, $\oper{R}{p}{n}_1$, is a restriction from the space of microscopic ensembles $D_\particlesub$ onto that of macroscopic ensembles $D_\noisysub$, i.e., $\oper{R}{p}{n}_1$ discards the information on the fast variable to obtain an ensemble $\mathcal{Z}$:
\begin{align*}
\oper{R}{p}{n}_1 : D_\particlesub \rightarrow D_\noisysub:
\mathcal{X} = \ensemble{X^p, Y^p, W^p}{p=0,1,\dots,P-1} \mapsto \mathcal{Z} = \ensemble{Z^p=X^p, W^p}{p=0,1,\dots,P-1}.
\end{align*}
The operator $\mathcal{R}_2$ measures the empirical distribution of an ensemble $\mathcal{Z}$ by binning into a histogram:
\begin{align*}
\oper{R}{n}{c}_2 : D_\noisysub \rightarrow D_\coarsesub:
\mathcal{Z} = \ensemble{Z^p, W^p}{p=0,1,\dots,P-1} \mapsto \vec{\rho}, 
\end{align*}
with 
\begin{align}
\vec\rho = (\rho_0, \rho_1, \dots, \rho_{I-1}), \qquad \text{with} \qquad \rho_i = \sum_{p \text{ s.t. } x_i < Z^p < x_{i+1}} W^p.
\label{eq:binning}
\end{align}
This procedure can be viewed as a statistical analogue of~\eqref{eq:integral}. 
Alternatively, one could consider a kernel density estimation (KDE) \citep{Scott15}, see Remark~\ref{remark:KDE} below. 
%The smoothing introduced by this this approach is moderate and certainly tolerable if the kernel width is chosen sufficiently small.

\paragraph{Matching.} 
Transferring coarse-level corrections to the microscopic level is much more involved, since there are many microscopic states consistent with a given macroscopic state.  We will define \emph{matching} operators that approach the problem as an inference problem: given a prior microscopic state $\mathcal{X}^*$, construct a posterior microscopic state $\overline{\mathcal{X}}$ that is consistent with the given macroscopic state $\vec\rho$. We denote matching as 
\[
\overline{\mathcal{X}}=\mathcal{M}(\vec\rho,\mathcal{X}^*).
\]
Practically speaking, we will choose $\mathcal{M}$ to minimise the perturbation with respect to the prior, in a sense that will be made precise in Section~\ref{sec:matching}. Some generic requirements on $\mathcal{M}$ will be given in Section~\ref{sec:couplings:requirements}. Specific choices for the matching operator will be given in Section~\ref{sec:matching}. One of the main objectives of this work is precisely to study the impact of this choice on parareal performance, a question will be discussed in Section~\ref{sec:results}. 

%\begin{remark}[On distributions and ensembles]
\medskip

In the above definition, the matching operator generates a posterior microscopic ensemble, given a prior ensemble. It is often more convenient to think about matching as generating posterior \emph{distribution} based on the prior \emph{distribution}. Then, the matching operator $\mathcal{M}$ can be viewed as a two-step procedure: matching of distributions, denoted by $\mathcal{M}_1$, followed by resampling,  denoted by $\mathcal{M}_2$, i.e., 
\begin{align}
\oper{M}{}{} \left( \vec\rho, {\mathcal{X}^*} \right) = \oper{M}{}{}_2 \left( \oper{M}{}{}_1 \left( \vec\rho, \mu^* \right), {\mathcal{X}^*} \right),
\label{eq:shorthand}
\end{align} 
where the distribution $\mu^*$ is the distribution that is sampled by the prior ensemble $\mathcal{X}^*$. The operator $\mathcal{M}_2$ then represents an adjustment of the ensemble $\mathcal{X}^*$ to make it consistent with the distribution that is generated via the matching $\oper{M}{}{}_1 \left( \vec\rho, \mu^* \right)$. Studying the consistency of matching operators is most convenient at the distribution level, see, e.g., \citet{DeSaZi15,LeSaZi17}. In this work, we only consider matching operators for which we can perform the matching directly on the ensemble.
%\end{remark}

\subsubsection{Generic iteration ($k>0$)} \label{sec:micmac_parareal_generic}

With the propagators and matching/restriction operators available, we are ready to describe a generic micro-macro parareal iteration, i.e. any iteration except for the first one.
Each parareal iteration starts from an approximation of the full microscopic state $\left({\mathcal{X}}_{k}^{n}\right)_{n=0}^N$ at iteration $k$ and at the time instances $t^n=n\Delta t$ ($0\le n \le N$), and corrects this microscopic state to obtain the approximation $\left({\mathcal{X}}_{k+1}^{n}\right)_{n=0}^N$ at iteration $k+1$. 

First, we compute, via restriction, the macroscopic approximation $\vec\rho_{k}^n = \oper{R}{}{}\mathcal{X}_{k}^{n}$, corresponding to the current microscopic approximation $\left({\mathcal{X}}_{k}^{n}\right)_{n=0}^N$.
After that, we perform $N$ parallel simulations using the fine propagator, using as the initial condition each current approximation $\mathcal{X}_{k}^{n}$, for $0 \le n \le N-1$. Simultaneously, we also perform $N$ parallel simulations using the coarse propagator, using as the initial condition the corresponding restricted state $\vec\rho_{k}^n$. We write this simulation phase as
\begin{align*}%\label{eq:mm_parallel}
\overline{\mathcal{X}}_{k}^{n+1} = \fineDt \left(\mathcal{X}_k^n\right), \qquad \text{and}\qquad \overline{\vec\rho}_k^{n+1}=\coarseDt \left(\vec\rho_{k}^n\right)=\coarseDt \left(\oper{R}{}{}\mathcal{X}_{k}^{n}\right)
\end{align*}
for any $0 \leq n \leq N-1$.

After these parallel computations, we have intermediate results $\overline{\mathcal{X}}_{k}^{n+1}$ at the fine level and $\overline{\vec\rho}_k^{n+1}$ at the coarse level, which are not necessarily consistent with each other. These intermediate results are then merged in a serial updating step, requiring communication between the processes. Given that only the macroscopic state is available at both the coarse and the fine levels, the corrections can only be performed on the macroscopic state. If we maintain the additive updating procedure of the classical parareal algorithm (as presented in Section~\ref{sec:parareal_traditional}), the micro-macro version of the serial update~\eqref{eq:sequential} can be written as
\begin{align}\label{eq:micmac_sequential_additive}
\vec\rho_{k+1}^{n+1} &= \coarseDt (\vec\rho_{k+1}^n) + \oper{R}{}{} \overline{\mathcal{X}}_{k}^{n+1} - \overline{\vec\rho}_k^{n+1}\\
 &= \coarseDt (\vec\rho_{k+1}^n) + \oper{R}{}{}\fineDt \left(\mathcal{X}_k^n\right)  - \coarseDt \left(\oper{R}{}{}\mathcal{X}_{k}^{n}\right) \nonumber
\end{align}
for any $0 \leq n \leq N-1$.

However, equation~\eqref{eq:micmac_sequential_additive} does not necessarily preserve positivity of the distributions. Therefore, we introduce a generalized operator $\jump{}$ that determines the macroscopic state of the next iterate, based on a sequential update and the use of the intermediate results $\overline{\mathcal{X}}_{k}^{n+1}$ and $\overline{\vec\rho}_k^{n+1}$, i.e.,
\begin{align}
\vec\rho_{k+1}^{n+1} &= \jump{}\left(
\coarseDt (\vec\rho_{k+1}^n), \,
\oper{R}{}{} \overline{\mathcal{X}}_{k}^{n+1}, \,
\overline{\vec\rho}_k^{n+1} \,
\right) \nonumber 
\\
&= \jump{}\left(
\coarseDt (\vec\rho_{k+1}^n), \,
\oper{R}{}{}\fineDt \left(\mathcal{X}_k^n\right) , \,
\coarseDt \left(\oper{R}{}{}\mathcal{X}_{k}^{n}\right) \,
\right)
\label{eq:mm_sequential}
\end{align}
for any $0 \leq n \leq N-1$, with the initial condition $\vec\rho_{k+1}^{0} = \vec\rho_k^{0}$.

\begin{remark}[On the arguments of $\mathcal{J}$]\label{rem:argJ}
The operator $\mathcal{J}$ computes the coarse state $\vec\rho_{k+1}^{n+1}$ at iteration $k+1$ and time $n+1$, based on three arguments: (i) the coarse propagation of the coarse state at iteration $k+1$ and time $n$, denoted as $\coarseDt (\vec\rho_{k+1}^n)$; (ii) the restriction of the fine intermediate state $\mathcal{R}\overline{\mathcal{X}}_{k}^{n+1}$; and (iii) the coarse intermediate state $\overline{\vec\rho}_k^{n+1}$. In Section~\ref{sec:jump_convergence} and Section~\ref{sec:jumps}, we will study the properties and choices of $\mathcal{J}$ in isolation, without coupling to the iteration number $k$ and the time instance $n$.  We then write the serial updating step as 
\begin{equation}\label{eq:jump-notation}
\vec\rho_{k+1}^{n+1}=\mathcal{J}(\vec\rho_{\textrm{c}}^{n+1},\vec\rho_{\textrm{f}}^{*},\vec\rho_{\textrm{c}}^{*}),
\end{equation}
i.e., we denote the coarse state that is to be corrected as $\vec\rho_{\textrm{c}}^{n+1}=\coarseDt (\vec\rho_{k+1}^n)$, and the two intermediate values that will be used for the correction as $\vec\rho_{\textrm{f}}^{*}=\mathcal{R}\overline{\mathcal{X}}_{k}^{n+1}=\oper{R}{}{}\fineDt \left(\mathcal{X}_k^n\right)$ and $\vec\rho_{\textrm{c}}^{*}=\overline{\vec\rho}_k^{n+1}=\coarseDt \left(\oper{R}{}{}\mathcal{X}_{k}^{n}\right)$.
\end{remark}

%Note that the coarse propagator result $\coarseDt\left(\vec\rho_k^n\right)$ does not have to be computed at this iteration, as it is was already computed in the previous iteration for the first argument of $\jump{}$. 
In Section~\ref{sec:jumps}, we describe in detail several options for choosing $\jump{}$. The effect of the choice of $\jump{}$ on the performance of the micro-macro parareal method will be discussed in Section~\ref{sec:results}, where we perform comprehensive numerical experiments.

At this point, we now only have a corrected \emph{macroscopic} state $\vec\rho_{k+1}^{n+1}$ given by~\eqref{eq:micmac_sequential_additive} or~\eqref{eq:mm_sequential}. To obtain a corresponding \emph{microscopic} state $\left({\mathcal{X}}_{k+1}^{n}\right)_{n=0}^N$, we perform a matching step, using the intermediate microscopic state $\overline{\mathcal{X}}_{k}^{n+1}$ as the prior,
\begin{align}
\mathcal{X}_{k+1}^{n+1} = \oper{M}{}{} \left( \vec\rho_{k+1}^{n+1}, \overline{\mathcal{X}}_k^{n+1} \right).
\label{eq:mm_match}
\end{align}
This step concludes one generic iteration.

\begin{remark}
%\label{remark:viewpoint}
The order of the steps within a single micro-macro parareal iteration depends on our viewpoint. We may view the algorithm in two ways: (i) as a way to accelerate fine simulations, as was done in the description above; or (ii) as a way of increasing the fidelity of coarse simulations. In the second viewpoint, one micro-macro parareal iteration iteratively corrects the \emph{macroscopic} state.  
There is no difference in the implementation of the algorithm in the alternative viewpoint, but it does change what is considered the last step of the current iteration and what is the first step of the next iteration. (In viewpoint (ii), the matching constitutes the start of a new iteration.)
\end{remark}

\subsubsection{Zeroth iteration ($k=0$)} \label{sec:micmac_parareal_zero}

In the traditional parareal algorithm of Section~\ref{sec:parareal_traditional}, the zeroth iteration is given by a simple serial pass of the coarse propagator over the entire time domain. In the micro-macro parareal algorithm, this step only gives us an approximation of the \emph{macroscopic} state at each time instance,
\begin{align*}
\vec\rho_{k=0}^{n+1} = \coarseDt (\vec\rho_{k=0}^n) \qquad \mbox{  with  } \qquad \vec\rho_{k=0}^0 = \oper{R}{}{} \mathcal{X}_0,
%\label{eq:mm_initial}
\end{align*}
in which $\mathcal{X}_0$ represents the (microscopic) initial condition. We cannot perform the matching step~\eqref{eq:mm_match} 
to obtain a corresponding microscopic state, because no prior microscopic state is available from a previous iteration at the corresponding moment in time. We resolve this issue by performing a matching with respect to the initial condition $\mathcal{X}_0$ in the zeroth iteration, i.e.,
\begin{align*}
\mathcal{X}_{k=0}^n = \oper{M}{}{} \left( \vec\rho_{k=0}^n, \mathcal{X}_0 \right), \qquad 1 \le n \le N.
%\label{eq:mm_lifting}
\end{align*}
Other choices are possible, such as matching with the global equilibrium solution \citep{LeLeSa13} or using some entropy maximization \citep{SaLeLe11}.

\subsection{Properties of the coupling operators}
\label{sec:couplings:requirements}
\citet{LeLeSa13} showed that convergence of a micro-macro parareal method for deterministic multiscale systems requires the matching to satisfy two generic properties: consistency and a projection property. In this section, we first discuss consistency (Section~\ref{sec:consist}) and the projection property (Section~\ref{sec:proj}) in the slow-fast SDE setting. We next comment on the consequences of these properties on the behaviour of the micro-macro parareal method in Section~\ref{sec:jump_convergence}.

\subsubsection{Consistency}\label{sec:consist}

The \emph{consistency} property \citep[Eq.~(3.3)]{LeLeSa13} requires that the microscopic state that results from the matching is indeed consistent with the imposed macroscopic state. With the above notation, consistency writes
\begin{align}
\vec\rho = \left( \oper{R}{p}{n} \circ \oper{M}{c}{p} \right) \left( \vec\rho, \mathcal{X} \right) \mbox{ for any } \mathcal{X} \in D_\particlesub \mbox{ and } \vec\rho \in D_\coarsesub,
\label{eq:strong_consistency}
\end{align}
i.e., the consistency property holds exactly for a given ensemble $\mathcal{X}$. In the (stochastic) Monte Carlo setting, we call~\eqref{eq:strong_consistency} \emph{strong consistency}. This \emph{deterministic} requirement may however be too strict and preclude the use of sampling procedures that have otherwise advantageous properties, see Section~\ref{sec:matching}. 

Therefore, we also formulate a \emph{weak consistency} property, in which we only require consistency in expectation, 
\begin{align}
\expbra{ \vec\rho} = \expbra{ \left( \oper{R}{p}{n} \circ \oper{M}{c}{p} \right) \left( \vec\rho, \mathcal{X} \right) } \mbox{ for any } \mathcal{X} \in D_\particlesub \mbox{ and } \vec\rho \in D_\coarsesub,
\label{eq:consistency}
\end{align}
where the expectation is taken with respect to the random process used in the sampling operator $\oper{M}{f}{p}_2$ introduced in~\eqref{eq:shorthand}.
Both strong and weak consistency can be realized with an appropriate sampling procedure $\mathcal{M}_2$ in the matching operator. How strong and weak consistency impact the micro-macro parareal algorithm will be discussed in Section~\ref{sec:matching}.
\begin{remark}
\label{remark:KDE}
The use of kernel density estimation, see e.g.~\cite{Scott15}, instead of a histogram approximation for taking the empirical distribution violates both consistency properties, leading us to not include it further in our exposition.
\end{remark}

\subsubsection{Projection}\label{sec:proj}

The projection property \citep[Eqs.~(3.6)--(3.7)]{LeLeSa13} requires that, if the prior used in the matching is already consistent with the imposed macroscopic state, the resulting matched microscopic state equals the prior. 
We can again require this property to hold exactly for a given ensemble $\mathcal{X}$, which we write as 
\begin{align}
\oper{M}{c}{p} \left( \oper{R}{p}{n} \mathcal{X}, \mathcal{X} \right) = \mathcal{X} \mbox{ for any } \mathcal{X} \in D_\particlesub.
\label{eq:projection}
\end{align}
%This property ensures that once a parareal iteration has converged for a certain time step, the ensuing fine scale simulation is initialized with an identical ensemble resulting in an identical simulation, if care is taken to reproduce the same Brownian paths. 
%In practice, of course, such simulations are redundant and are not recomputed.
We call~\eqref{eq:projection} the \emph{strong projection property}. A weak version of the projection property would look at only the matching operator (and not the sampling) and state
\begin{align*}
\left( \oper{M}{c}{p}_1 \right) \left( \vec\rho, \mu \right) = \mu \mbox{ for any } \mu \in D_\finesub \mbox{ with } \vec\rho = \int \mu \dd y.
\end{align*}
%This, however, negates the strong convergence to a single realization of the fine propagator and consequently disrupts the basis for omitting redundant simulations. 
%To maintain strong convergence, we demand the coupling strategies respect~\eqref{eq:projection}.
In the next Section~\ref{sec:jump_convergence}, we show that the strong formulation~\eqref{eq:projection} has advantageous properties. 

\subsubsection{Exactness of the multiscale parareal algorithm}
\label{sec:jump_convergence}

The micro-macro parareal algorithm outlined in this section satisfies a similar exactness property as the traditional algorithm of Section~\ref{sec:parareal_traditional}, provided the matching operator satisfies the strong projection property~\eqref{eq:projection} and the iterator $\jump{}$ (defined in~\eqref{eq:mm_sequential}) satisfies a strong consistency property,
\begin{align}
\mathcal{J}(\vec\rho_{\textrm{c}}^{n+1},\vec\rho_{\textrm{f}}^{*},\vec\rho_{\textrm{c}}^{*})=
\vec\rho_{\textrm{f}}^{*} \quad \mbox{ whenever } \quad \vec\rho_{\textrm{c}}^{n+1} \equiv \vec\rho_{\textrm{c}}^{*},
\label{eq:iterator_property}
\end{align}
in which we have used the notation of~\eqref{eq:jump-notation}. Then, the $k$-th iteration of the micro-macro parareal algorithm is exact for all times $n\Delta t$ with $n \leq k$, i.e., the microscopic state satisfies $\mathcal{X}_k^n = \left(\fineDt\right)^n \mathcal{X}_0$ for all $n \leq k$. 
 
We first demonstrate this exactness for $n=k=1$ and then proceed by induction. The iterator~\eqref{eq:mm_sequential} for $n=1$ in the first iteration $k=1$ reads
\begin{align*}
\vec\rho_{k=1}^{n=1} = \jump{}\left(
\coarseDt (\vec\rho_{k=0}^{n=0}), \,
\oper{R}{}{} \overline{\mathcal{X}}_{k=0}^{n=1}, \,
\coarseDt (\vec\rho_{k=1}^{n=0}) \,
%\vec\rho_k^{n+1}
\right),
\quad \mbox{with} \quad
\vec\rho_{k=1}^{n=0} = \vec\rho_{k=0}^{n=0} = \oper{R}{}{} \mathcal{X}_0.
%\label{eq:mm_sequential_first}
\end{align*}
As the initial condition for both coarse propagators ($\vec\rho_{k=0}^{n=0}$ and $\vec\rho_{k=1}^{n=0}$) are identical, so are the results $\coarseDt (\vec\rho_{k=0}^{n=0})$ and $\coarseDt (\vec\rho_{k=1}^{n=0})$.  
Using~\eqref{eq:iterator_property}, we therefore have $\vec\rho_1^1=\oper{R}{}{} \overline{\mathcal{X}}_{k=0}^{n=1}$. 
The strong projection property~\eqref{eq:projection} then guarantees that the result of the matching operator $\mathcal{X}_1^1$ is equal, in a strong sense, to the result of the microscopic propagation $\overline{\mathcal{X}}_{k=0}^{n=1} = \fineDt \mathcal{X}_0$. 
The exact result at $t=\Delta t$ after the first iteration is then used as a starting point to show that the solution after a second iteration is exact at $t=2\Delta t$, and so on. 

This exactness property also provides some computational gain \citep{Aubanel11}, since those fine simulations with $n \leq k$ need not be computed. 
For a simulation using $K$ iterations, $\frac{K(K-1)}{2}$ fine scale simulations need not be performed. 
However, since it is ultimately desirable that $K \ll N$, this gain is usually marginal. 

%Having discussed how the parareal algorithm can be applied to probabilistic multiscale problems, we now delve into some more detail on two of the components. 
%The matching and sampling, combined in practice in a single procedure, are discussed in more detail in Section~\ref{sec:matching}. 
%In Section~\ref{sec:jumps} we present more details on how the iterative step of~\eqref{eq:mm_sequential} is performed. 
%As the effect of the different implementations can only discussed once the procedures are fully developed, numerical results are postponed until Section~\ref{sec:results}.

\section{Iterative steps}
\label{sec:jumps}

The operator $\jump{}(\vec\rho_{\textrm{c}}^{n+1}, \vec\rho_{\textrm{f}}^{*}, \vec\rho_{\textrm{c}}^{*})$ in~\eqref{eq:mm_sequential} determines the macroscopic state of the next iteration. (See Remark~\ref{rem:argJ} for the meaning of the three arguments of $\mathcal{J}$.)
This operator is the multiscale equivalent of~\eqref{eq:sequential} in the traditional algorithm. 
The goal is to use the discrepancy between the result of the coarse propagator $\vec\rho_{\textrm{c}}^{*} = \coarseDt \left(\oper{R}{}{}\mathcal{X}_{k}^{n}\right)$ and the restriction of the result of the fine propagator $\vec\rho_{\textrm{f}}^{*} = \oper{R}{}{} \fineDt \mathcal{X}^n_k$ to update the current iteration of the coarse propagator $\vec\rho_{\textrm{c}}^{n+1} = \coarseDt \vec\rho^n_{k+1}$.
The iterators should respect positivity and unit mass of the distributions as well as the property~\eqref{eq:iterator_property}, as explained in Section~\ref{sec:parareal_micmac}.
Without positivity and unit mass, the matching operators of Section~\ref{sec:matching} are not effective%
\footnote{
The reweighting sampling strategy of Section~\ref{sec:matching:conditional} would assign a negative weight to certain particles if the macroscopic distribution were negative. While this is not necessarily problematic, we still choose to enforce positivity of the distribution produced by the iterative step.
}.
The third requirement (i.e. property~\eqref{eq:iterator_property}) is necessary to guarantee that when the iteration has converged, the iterator returns the macroscopic restriction of the fine propagator result, $\vec\rho_{\textrm{f}}^{*}$, as detailed in Section~\ref{sec:jump_convergence}. In this section, we detail four choices for the iterator and demonstrate their behaviour for three test cases. Numerical results are presented in Section~\ref{sec:results}. 

\subsection{Additive iterator}
%\label{sec:jumps:add}

A direct analogue to the traditional iteration in~\eqref{eq:sequential} is given by 
\begin{align*}
\jump{add}\left( \vec\rho_{\textrm{c}}^{n+1}, \vec\rho_{\textrm{f}}^{*}, \vec\rho_{\textrm{c}}^{*} \right) = \vec\rho_{\textrm{c}}^{n+1} + \vec\rho_{\textrm{f}}^{*} - \vec\rho_{\textrm{c}}^{*}.
%\label{eq:jump_add}
\end{align*}
In the numerical implementation of this iterator $\jump{add}$, we tackle lack of positivity by setting all negative values of the result to zero, and then rescaling to satisfy unit mass. 

%\begin{figure}[hpt]
%\begin{center}
%\drawjumppdfs{add}
%\caption{
%The result of the additive iterator for three different test cases. 
%In all three panels the colours are as follow: $\rho$ is in solid dark blue, $\overline\rho$ in solid orange, $\hat\rho$ in dashed light blue and the result $\tilde\rho$ is in dashed dark red.
%The left panel tests the situation where the two coarse propagator results are equal, $\vec\rho_{\textrm{c}}^{n+1} \equiv \vec\rho_{\textrm{c}}^{*}$, where we would expect to see $\tilde{\vec\rho} \equiv \vec{\overline\rho}$.
%In the centre panel $\hat\rho(x) = \rho(x-1.8)$, corresponding to a pure advection. 
%In the right panel $\hat\rho(x) = \rho(\frac{5}{3}x)$, corresponding to a sharpening distribution. 
%}
%\label{fig:add_jump}
%\end{center}
%\end{figure}
%Figure~\ref{fig:add_jump} displays the result of an iteration step $\tilde{\vec\rho} = \jump{add}\left( \vec\rho_{\textrm{c}}^{n+1}, \vec\rho_{\textrm{f}}^{*}, \vec\rho \right)$ for three different cases. 
%The left panel demonstrates that when the previous and current coarse propagator results ($\hat\rho$ and $\rho$) are equal, the result of the iteration is to return the restriction of the fine propagator result $\overline\rho$.
%For the centre and right panels there is not such a clearly defined correct result. 
%Nevertheless, it is interesting to note that on the centre panel, where the only difference between $\hat\rho$ and $\rho$ is an advection, the result shows bimodality. 
%Whether this is detrimental to results is tested in Section~\ref{sec:results}. 

\subsection{Multiplicative iterator} 
%\label{sec:jumps:mult}

The multiplicative iterator is given by
\begin{align*}
\jump{mult}\left( \vec\rho_{\textrm{c}}^{n+1}, \vec\rho_{\textrm{f}}^{*}, \vec\rho_{\textrm{c}}^{*} \right) = \vec\rho_{\textrm{c}}^{n+1} \ \ \frac{\vec\rho_{\textrm{f}}^{*}}{\vec\rho_{\textrm{c}}^{*}}.
%\label{eq:jump_mult}
\end{align*}
This iterator naturally maintains positivity, but not unit mass. 
Also, the function is not well defined if $\vec\rho_{\textrm{f}}^{*}$ is not supported by $\vec\rho_{\textrm{c}}^{*}$. 
This is remedied by setting the result to zero in those areas, and rescaling the solution accordingly.

%\begin{figure}[hpt]
%\begin{center}
%\drawjumppdfs{mult}
%\caption{
%The result of the multiplicative iterator for the same three test cases as in Figure~\ref{fig:add_jump}, using the same colour scheme.
%}
%\label{fig:mult_jump}
%\end{center}
%\end{figure}
%Figure~\ref{fig:mult_jump} displays the result of an iteration step $\tilde{\vec\rho} = \jump{mult}\left( \vec\rho_{\textrm{c}}^{n+1}, \vec\rho_{\textrm{f}}^{*}, \vec\rho_{\textrm{c}}^{*} \right)$ for three different cases. 
%The left panel demonstrates that when the previous and current coarse propagator results ($\hat\rho$ and $\rho$) are equal, the result of the iteration is to return the restriction of the fine propagator result $\overline\rho$.
%For the centre and right panels there is not such a clearly defined correct result. 
%Nevertheless, it is interesting to note that on the centre panel, there is no bimodality in the result, as opposed to the additive iterator of Section~\ref{sec:jumps:add}. 
%The numerical results in Section~\ref{sec:results} show how this impacts the quality of the parareal algorithm. 

\subsection{Rotation iterator}
%\label{sec:jumps:rotate}

The third approach uses rotations of vectors in the Hilbert space $L^2(\IR)$, and follows the idea suggested in~\citep[Sec.~3]{MaTu02} when applying the parareal algorithm to the Schroedinger equation (in that setting, preserving the $L^2$ norm of the solution is critical). This requires here embedding the probability distributions $\vec\rho_{\textrm{c}}^{n+1}$, $\vec\rho_{\textrm{f}}^{*}$ and $\vec\rho_{\textrm{c}}^{*}$, which all belong to $L^1(\IR)$, in $L^2(\IR)$. Due to the non-negativity of probability distributions, this embedding is possible by taking the square root: $\rho \mapsto \sqrt{\rho}$. This embedding puts probability density functions onto the unit sphere of $L^2(\IR)$, where we may define a rotation (in the plane spanned by $\sqrt{\vec\rho_{\textrm{c}}^{*}}$ and $\sqrt{\vec\rho_{\textrm{c}}^{n+1}}$) that brings $\sqrt{\vec\rho_{\textrm{c}}^{*}}$ to $\sqrt{\vec\rho_{\textrm{c}}^{n+1}}$. By applying this rotation to $\sqrt{\vec\rho_{\textrm{f}}^{*}}$, we find our result in terms of a vector on the $L^2(\IR)$ unit sphere. Squaring this result to map back to a norm one, non-negative function in $L^1(\IR)$ yields
\begin{align*}
\jump{rot}\left( \vec\rho_{\textrm{c}}^{n+1}, \vec\rho_{\textrm{f}}^{*}, \vec\rho_{\textrm{c}}^{*} \right) = \left( 
\frac{u+2uw-v}{1+w} \sqrt{ \vec\rho_{\textrm{c}}^{n+1} } + 
\sqrt{ \vec\rho_{\textrm{f}}^{*} } - 
\frac{u+v}{1+w} \sqrt{ \vec\rho_{\textrm{c}}^{*} }
\right)^2,
%\label{eq:jump_rot}
\end{align*}
where
\begin{align*}
u &= \int \sqrt{ \vec\rho_{\textrm{f}}^{*} } \ \sqrt{ \vec\rho_{\textrm{c}}^{*} } = \sum_{i=0}^{I-1} \sqrt{ \vec\rho_{\textrm{f},i}^{*} \ \vec\rho_{\textrm{c},i}^{*}},  \\
v &= \int \sqrt{ \vec\rho_{\textrm{f}}^{*} } \ \sqrt{ \vec\rho_{\textrm{c}}^{n+1} } = \sum_{i=0}^{I-1} \sqrt{ \vec\rho_{\textrm{f},i}^{*} \ \vec\rho_{\textrm{c},i}^{n+1}},\\
w &= \int \sqrt{ \vec\rho_{\textrm{c}}^{*}} \ \sqrt{ \vec\rho_{\textrm{c}}^{n+1} } = \sum_{i=0}^{I-1} \sqrt{ \vec\rho_{\textrm{c},i}^{*} \ \vec\rho_{\textrm{c},i}^{n+1} }, 
\end{align*}
in which the sum is taken over the values of the densities in the $I$ bins of the histogram, see equation~\eqref{eq:piecewise_constant_pdf}. 
This approach naturally maintains both positivity and unit mass. 
We also readily check that, if $\vec\rho_{\textrm{c}}^{*} \equiv \vec\rho_{\textrm{c}}^{n+1}$, then the rotation reduces to the identity operator ($w = 1$ and $u = v$) and $\jump{rot}\left( \vec\rho_{\textrm{c}}^{*}, \vec\rho_{\textrm{f}}^{*}, \vec\rho_{\textrm{c}}^{*} \right) = \vec\rho_{\textrm{f}}^{*}$, as required in~\eqref{eq:iterator_property}.

%\begin{figure}[hpt]
%\begin{center}
%\drawjumppdfs{rot}
%\caption{
%The result of the rotation iterator for the same three test cases as in Figure~\ref{fig:add_jump}, using the same colour scheme.
%}
%\label{fig:rotate_jump}
%\end{center}
%\end{figure}
%Figure~\ref{fig:rotate_jump} displays the result of an iteration step $\tilde{\vec\rho} = \jump{rot}\left( \vec\rho_{\textrm{c}}^{n+1}, \overline{\vec\rho}, \vec\rho \right)$ for three different cases. 
%The left and right hand panels show that this matching behaves similar to the addition and multiplication iterators of the previous Sections. 
%In the centre panel, there is a slight bimodality, but much less pronounced than in the addition iterator. 
%Note that this iterator natively satisfies unit mass (as does the addition iterator) and positivity (as does the multiplication iterator), so perhaps the intermediate result of the advection test case is unsurprising. 

\subsection{Quantile iterator}
%\label{sec:jumps:cdf}

The fourth iterator is a simple addition and subtraction, but applied to the quantile function of the distribution $\rho$ rather than the distribution itself. A similar technique is used by \citet{Gear01}. The quantile function $q(p), p \in \left[ 0,1 \right]$, is given by the inverse of the cumulative density function (cdf) of $\rho$, denoted $c$ and defined by
\[
c(x) = \int_{-\infty}^x \rho, \qquad q(p) = c^{-1}(p) \quad \text{for any} \quad p\in[0,1].
\] 
Note that, being based on the cumulative distribution, this approach cannot readily be generalised to the case of a \emph{slow} variable in dimension larger than one.

%As explained for the matching in Section~\ref{sec:matching:conditional} we may assume that, in the absence of statistical error, all distributions used have infinite support. 

We denote the three quantile functions, corresponding to the three densities $\vec\rho_{\textrm{c}}^{n+1}$, $\vec\rho_{\textrm{f}}^{*}$ and $\vec\rho_{\textrm{c}}^{*}$, as respectively $q_{\textrm{c}}^{n+1}$, $q_{\textrm{f}}^{*}$ and  $q_{\textrm{c}}^{*}$.  The quantile iterator is then expressed by first computing the quantile function
\begin{align}\label{eq:quantile_reconstruction}
\forall p \in [0,1], \quad q_{k+1}^{n+1}(p) = q_{\textrm{c}}^{n+1}(p) + q_{\textrm{f}}^{*}(p) - q_{\textrm{c}}^{*}(p)
\end{align}
and then taking the derivative of the inverse of $q_{k+1}^{n+1}$ to yield $\vec\rho^{n+1}_{k+1}$.

In the discretised setting, we start from the vector $\vec\rho$ (see equation~\eqref{eq:piecewise_constant_pdf}) that represents a piecewise constant probability density function to define a piecewise linear cumulative density function
\begin{align}\label{eq:discr_cdf}
c(x) = \sum_{i=0}^I \bs{1}_{(x_i,x_{i+1}]}(x) \left( \frac{x - x_i}{x_{i+1} - x_i} \rho_i + \sum_{j=0}^{i-1} \rho_{j} \right),
\end{align}
which we use for all densities involved. 
For the empirical distribution that follows from the restriction of $\overline{\mathcal{X}}$, there also exists a direct empirical cdf, given by $\overline{c}_\mathrm{emp}(x) = \sum_{p \text{ s.t. } X^p < x} W^p$. 
This empirical cdf is, however, more expensive to compute than~\eqref{eq:discr_cdf} and is not used explicitly. Note also that the value of the piecewise linear cdf ${c}$ evaluated at the bin edges $x_i$ is exact, i.e., for any $i$, we have
\begin{align*}
c(x_i) = \int_{-\infty}^{x_i} \rho(x) \dd x = \sum_{j=0}^{i-1} \rho_j = \sum_{j=0}^{i-1} \ \sum_{p \text{ s.t. } x_j < X^p < x_{j+1}} W^p = \sum_{p \text{ s.t. } X^p < x_i} W^p = \overline{c}_\mathrm{emp}(x_i).
%\label{eq:exactness_cdf}
\end{align*}

To implement the quantile reconstruction~\eqref{eq:quantile_reconstruction}, we start from the evaluated empirical cdf at the bin edges, which we write as tuples $(x_i,c(x_i))=(x_i,p_i)=(q_i,p_i)$. Inverting the cdf to obtain the quantile function amounts to writing these tuples as $(p_i,q_i)$. Since each of the three quantile functions involved is evaluated on the same set of histogram bins, they are not available on the same set of points $p_i$. We choose to use the \emph{union} of these values, i.e., we take the set of the $\tilde{N}$ \emph{distinct} values $p_j$ that any of the involved cdf's takes at the bin edges $\{ p_j \}_{j=0,\dots,\tilde{N}}$, ordered such that $p_{j+1} > p_j$.
The quantile functions $\hat{q}$, $\overline{q}$ and $q$, corresponding to $\vec\rho_{\textrm{c}}^{n+1}$, $\vec\rho_{\textrm{f}}^{*}$ and $\vec\rho_{\textrm{c}}^{*}$ respectively, are evaluated at these values $p_j$ and then added and subtracted to find
\begin{align*}
\tilde{q}_j = \hat{q}(p_j) + \overline{q}(p_j) - q(p_j) \quad \mbox{for} \quad j=0,\dots,\tilde{N}.
\end{align*}
The resulting quantile function is defined by linearly interpolating between the points $\left( \tilde{q}_j, p_j \right)$. 
By finding the intersections of these lines with the bin edges, we retrieve the values of the resulting cdf at the bin edges $\tilde{c}_i$. 
The final result of the quantile iterator is the histogram $\tilde{\rho}_i$, readily found as the difference $\tilde{c}_{i+1} - \tilde{c}_i$.

Two practical issues need to be dealt with in the implementation of the quantile iterator. First, because the approximated probability density functions $\vec\rho_{\textrm{c}}^{n+1}$, $\vec\rho_{\textrm{f}}^{*}$ and $\vec\rho_{\textrm{c}}^{*}$ do not necessarily have infinite support, it is possible that the cumulative density functions are not invertible. 
In these cases, the quantile function is defined as follows:
\begin{itemize}
\item If the support of the pdf is bounded to the \emph{left}, we define $q(0) = \max \left\lbrace x \text{ s.t. } c(x) = 0 \right\rbrace$. 
If it is unbounded, $q(0)$ is not defined.
\item If the support of the pdf is bounded to the \emph{right}, we define $q(1) = \min \left\lbrace x \text{ s.t. } c(x) = 1 \right\rbrace$. 
If it is unbounded, $q(1)$ is not defined.
\item For interior values $p \in (0,1)$ for which there is a range of values $x$ such that $c(x) = p$, we choose the average between the extremes, setting $q(p) = \half \min \left\lbrace x \text{ s.t. } c(x) = p \right\rbrace + \half \max \left\lbrace x \text{ s.t. } c(x) = p \right\rbrace$. 
\end{itemize}
%In fact, this last definition also holds for other interior values $p \in (0,1)$ where the cdf \emph{is} invertible.
Second, after adding and subtracting quantile functions, a lack of monotonicity may result in the obtained quantile function, i.e.~ $\tilde{q}_{j+1} < \tilde{q}_j$.
This occurs when the discrepancy between $\rho$ and $\overline\rho$ is such that it steepens an already steep cdf $\hat{c}(x)$. %, see the right hand panel of Figure~\ref{fig:cdf_jump} for a demonstration.
This is resolved by setting such quantiles $\tilde{q}_{j+1}$ and $\tilde{q}_j$ to their average.

\section{Matching and sampling}
\label{sec:matching}

This section details a number of options for the matching ($\oper{M}{c}{f}_1$) and sampling ($\oper{M}{f}{p}_2$) operators introduced in~\eqref{eq:shorthand}. 
In practice, these operators are combined, and the intermediate microscopic distribution $\mu$ is not explicitly computed. 
However, we still separately present the operators here, as doing so sheds light on the different implementations. 

We consider two different matching operators. Both produce a posterior microscopic distribution $\mu$ consistent with the target macroscopic distribution $\vec\rho$. The first is constructed such that the posterior \emph{marginal} distribution of the fast variable $y$ is identical to the prior, see Section~\ref{sec:matching:marginal}. The second leaves the \emph{conditional} distribution in $y$ for each given $x$ unchanged, see Section~\ref{sec:matching:conditional}.
%The matchings use a prior distribution $\overline\mu$, which has $\overline\rho(x) = \int_{-\infty}^\infty \mu(x,y) \dd y$ as its marginal in the slow variable.
Concurrently, we present the sampling procedures used for each of the matchings.

\subsection{Matching preserving the fast marginal}
\label{sec:matching:marginal}

This matching procedure is based on \emph{inverse transformation sampling} \citep{Devroye86}, which is commonly used for sampling arbitrary distributions. 
By applying the procedure to adjust the slow component of the prior ensemble, we leave the fast component untouched.
This ensures that the marginal distribution of the posterior remains the same as that of the prior. 
Given samples $\overline{X}^p$ from a distribution with pdf $\overline\rho(x)$ and cdf $\overline{c}(x)$, the procedure produces samples of a distribution $\rho$ and cdf $c$ via the transformation
\begin{align}
X^p = c^{-1}\left( \overline{c}(\overline{X}^p ) \right), \quad Y^p = \overline{Y}^p.
\label{eq:transformation}
\end{align}
This transformation is the solution to an optimal transport problem in dimension one. It cannot readily be generalized to the case of a \emph{slow} variable in larger dimension.
It works by first transforming the arbitrary distribution of $\overline{X}^p$ to a uniform distribution via the cdf $\overline{c}$ and then transforming these uniformly distributed samples to the distribution with pdf $\rho$ via the inverse cdf $c^{-1}$. 
Corresponding to this transformation of the ensemble is a distribution sampled by the new ensemble.
By applying a coordinate transformation we find its form
\begin{align*}
\mu(x,y) = \overline\mu\left(\overline{c}^{-1}\left(c(x)\right), y\right) \ \frac{\rho(x)}{\overline{\rho}\left(\overline{c}^{-1}\left(c(x)\right)\right)}
\end{align*}
in terms of the prior distribution $\overline\mu$ sampled by the prior ensemble $\overline{\mathcal{X}}$. 
It is readily verified that the marginal of the slow variable $\rho(x) = \int \mu(x,y) \dd y$ is as desired and that the marginal of the fast variable $\phi(y) = \int \mu(x,y) \dd x$ is equal to the prior marginal $\overline\phi(y) = \int \overline\mu(x,y) \dd x$.

The sampling procedure is evident for this matching, it has already been stated in~\eqref{eq:transformation}. 
The empirical cumulative density function $\overline{c}$ of the prior ensemble is given by a piecewise constant function
\begin{align*}
\overline{c}(x) = \sum_{p \text{ s.t. } \overline{X}^p < x} \overline{W}^p.
%\label{eq:exact_cdf}
\end{align*}
Note that in this matching strategy, the weights $W^p$ remain unchanged throughout the iterations. Nevertheless, if the initial condition has unequal weights (importance sampling), the $W^p$ are not equal. 

To avoid the computational cost of sorting the particles, which is required for computing this cumulative density, we use instead the integrated empirical probability density function from~\eqref{eq:binning} and~\eqref{eq:piecewise_constant_pdf} as an approximation. 
Evaluating the piecewise linear approximated cdf at the position of each particle $X^p$ results in
\begin{align*}
\overline{c}^p = \sum_{j=0}^{i-1} \overline\rho_j + \frac{\overline{X}^p - x_i}{x_{i+1} - x_i}\rho_i \approx \overline{c}(\overline{X}^p),
%\label{eq:modulation_quantiles}
\end{align*}
%\todo{Check this formula (bars?)}
where $i$ is the bin number such that $x_i < \overline{X}^p < x_{i+1}$. 
The inversion of $c$ is done using a similar piecewise linear cumulative density function following from the posterior macroscopic distribution $\rho$. 
This inversion takes the form
\begin{align*}
X^p = x_i + \frac{\overline{c}^p - c_i}{\overline\rho_i} \left(x_{i+1} - x_i\right), 
\end{align*}
where $i$ is the bin number such that $c_i < \overline{c}^p < c_{i+1}$ with $c_i = \sum_{j=0}^{i-1} \rho_j$. 
Note that the definition of the sampling does not include the posterior or prior microscopic distributions, but only the macroscopic distributions $\overline{\rho}$ and $\rho$. 
The performance of a parareal implementation with this matching and sampling procedure is compared against others in Section~\ref{sec:results}.

\subsection{Matching preserving the conditional distribution}
\label{sec:matching:conditional}

The conditional distribution of the fast variable is readily preserved if the posterior probability distribution $\mu$ resulting from the matching has the form
\begin{align*}
\mu(x,y) = \overline\mu(x,y) f(x).
\end{align*}
Of course, the function $f$ must be such that $\mu$ remains a valid distribution, i.e.~non-negative and with unit mass. 
The densities involved in solving the system of SDEs~\eqref{eq:SDEx}--\eqref{eq:SDEy} and in the Fokker-Planck equation~\eqref{eq:FPz} all have infinite support. 
Even if the initial conditions have bounded support, after an arbitrarily small time the solutions will have infinite support. 
Therefore we may assume that the posterior is supported by the prior, namely $\rho \ll \overline\rho$, and we can take $\displaystyle f(x) = \frac{\rho(x)}{\overline\rho(x)}$. 
By integrating $\mu$ in $y$ we find
\begin{align*}
\int \mu(x,y) \dd y = \frac{\rho(x)}{\overline\rho(x)} \int \overline\mu(x,y) \dd y = \frac{\rho(x)}{\overline\rho(x)} \ \overline\rho(x) = \rho(x),
\end{align*}
demonstrating that the posterior marginal in the slow variable is as required and also that the posterior has unit mass. 
Positivity is satisfied via the positivity of both prior and posterior marginals. 

We consider two sampling procedures while implementing this matching. The first is a simple reweighting of the prior ensemble $\overline{\mathcal{X}}$, that is
\begin{align*}
W^p = \overline{W}^p \ \frac{\rho_i}{\overline\rho_i},
\end{align*}
where $i$ is the bin number such that $x_i < \overline{X}^p < x_{i+1}$. 
In cases where the \emph{empirical} prior $\overline{\vec\rho}$ does not support the posterior $\vec\rho$, i.e.~when there exist histogram bins $i$ for which $\rho_i \neq \overline\rho_i = 0$, this reweighting would lead to a loss of mass. 
This loss of mass is remedied by moving the mass from bins that have no particles to the nearest bins that do, half of it to the left, half to the right.
This occurs mostly during the first iterations, as the solution converges, and only at the very tails of the distribution. 
%This can be remedied in two ways: 
%\begin{itemize}
%\item Scale the supported elements of $\vec\rho$ to recover the mass of the unsupported elements. This divides the mass over all ensemble members in proportion to their likelihood. 
%\item Divide the mass of the unsupported elements in two and add each half to the nearest supported element on the left and right. If for some unsupported element there is no supported element to the left or right, add all the mass to the side where there is a supported element.
%\end{itemize}
%The former option introduces a bias, as elements that are likely to be well-supported are usually scaled to carry more weight, whereas elements that are often not supported carry too little weight, on average. 
%The latter also has a bias, but this bias only occurs near the extemes of the prior distribution, where it is likely that the mass cannot be spread evenly to both sides. 
%See Section~\ref{sec:matching:consistency:conditional} for an example of this bias. 

In the second sampling strategy, the reweighting is followed by a resampling step. 
We use a \emph{stratified resampling} technique \citep{DoCa05, HoScGu06, MaJoLe07}, whereby each particle $p$ from the prior ensemble $\overline{\mathcal{X}}$ occurs in the posterior $n_p$ times, with
\begin{align*}
n_p &= \#\left\lbrace u_j \ \text{s.t.} \ \sum_{q=0}^{p-1} W^q < u_j < \sum_{q=0}^{p} W^q \right\rbrace,
%n_p &= \#\left\lbrace u_j \,:\, 0 < u_j - \sum_{q=0}^{p-1} W^q < W^p \right\rbrace,
\end{align*}
where the $u_j$ are random numbers following
\begin{align}
u_j = \frac{j + \tilde{u}_j}{P}, \qquad \tilde{u}_j \sim U[0,1) \mbox{ i.i.d}.
\label{eq:stratified}
\end{align}
The stratified form of this resampling is desirable as it means small perturbations to the weights and thus leads to a small perturbation in the posterior ensemble. 
In the case of the weights remaining equal, the exact ensemble is reproduced, guaranteeing the strong projection property~\eqref{eq:projection}. 
These properties -- small perturbations to the ensembles for small changes in the macroscopic state and its limit case of an identical macroscopic state -- are particularly useful when interested in strong convergence. 
Note that again the definition of the sampling does not include the posterior or prior microscopic distributions, only the macroscopic distributions $\overline{\rho}$ and $\rho$. 

In the simulations, we will also consider a deterministic alternative to~\eqref{eq:stratified}, namely the same algorithm but with the deterministic choice $\tilde{u}_j \equiv \half$. 
This is a computationally cheaper approach that numerically appears to lead to equivalent results.

\section{Numerical results}
\label{sec:results}

In this section, we present numerical results for the parareal algorithm. Throughout the experiments, the coarse simulation uses the finite volume scheme described in Section~\ref{sec:mac_fv} on a uniform spatial grid of the space interval $[0;600]$ with the spacing $\Delta x = 1.0$, and the time step $\delta t_{\coarsesub} = 10^{-1}$. The microscopic ensemble is simulated using the Euler-Maruyama scheme with the time step $\delta t_{\finesub} = 10^{-3}$. Initial conditions are given by identical lognormal distributions for $X_0$ and $Y_0$, with mean $\mu_{X_0} = \mu_{Y_0} = 100$ and standard deviation $\sigma_{X_0} = \sigma_{Y_0} = 20$.
The ensemble size $P$ is specified for each of the test cases. 
We take care to use, in each of the time subintervals of the parareal simulation, the same random numbers in subsequent parareal iterations, so that the solution converges to a single realization of the $P$ Brownian processes.
%\footnote{The simulations with different time intervals do not represent identical Brownian processes.}
The value of the time-scale separation parameter $\epsilon$ is $2.5 \times 10^{-2}$ and the parameter $\tau_x$ is set to $25$.
These and other default parameters are collected in Table~\ref{tab:default_parameters} for later reference. 

\begin{table}[hpbt]
\caption{Default parameters used in the simulations}
\label{tab:default_parameters}
\begin{center}
\begin{tabular}{lcD{.}{.}{1.6}cD{.}{.}{1.6}} % dcolum package
%\multicolumn{2}{|c|}{Integers} &
%\multicolumn{2}{|c|}{}
Name & Symbol & \multicolumn{1}{c}{Value} & Reaction rate & \multicolumn{1}{c}{Value}
\\
\hline
Slow time scale  & $\tau_x$ & 2.5 \times 10^1 & $k_0$ & 6.0 \times 10^1
\\
Time-scale separation & $\epsilon$  & 2.5 \times 10^{-2} & $k_1$ & 1.0 
\\
Grid spacing & $\Delta x$ & 1.0 & $k_2$ & 4.8 \times 10^{-3}
\\
Coarse time step & $\delta t_{\coarsesub}$ & 1.0 \times 10^{-1} & $k_3$ & 0.666 \times 10^{-6}
\\
Fine time step & $\delta t_{\finesub}$ & 1.0 \times 10^{-3} & $k_4$ & 0.833
\\
Reporting time step & $\Delta t$ & 2.0 & $k_5$ & 1.0
\\
Number of intervals & $N$ & \multicolumn{1}{r}{20}  & % multicolumn to override the d-column decimal point placement for the integer values 
&
\\
Number of iterations & $K$ & \multicolumn{1}{r}{20} & &
\\
Ensemble size & $P$ & 1.0 \times 10^5 & &
\\
\hline
\end{tabular}
\end{center}
\end{table}

We perform a number of experiments to demonstrate the effectiveness of the algorithm in the current setting. 
In each case, we use all possible combinations of iterator strategies from Section~\ref{sec:jumps} and matching and sampling operators from Section~\ref{sec:matching}. 

\subsection{Varying the time scale separation}
\label{sec:results:epsilon}

In this first numerical experiment, we vary the time scale separation with values of $\epsilon$ ranging from $2.5 \times 10^{-2}$ down to $2.5 \times 10^{-3}$, leaving all other parameters as in Table~\ref{tab:default_parameters}. 
This changes the regime -- as explained in Section~\ref{sec:mic} -- of the dynamics between Regime 2 for modest time scale separation (e.g. $\epsilon = 2.5 \times 10^{-2}$) and Regime 1 for more extreme values. 
In the true limit of infinite time scale separation, the coarse model is accurate and the parareal algorithm is superfluous. This is the case when the dynamics is `safely' in Regime 1. 
However, for the values of $\epsilon$ that we consider here, we expect the coarse model to not be sufficiently accurate, which implies that the parareal algorithm is useful.
To test the performance of the algorithm in this setting, we consider a simulation over the time interval $[0,T]$ with $T=40$, with the initial conditions defined as above: mean $\mu_{X_0} = \mu_{Y_0} = 100$ and standard deviation $\sigma_{X_0} = \sigma_{Y_0} = 20$. 

\begin{figure}[hpt]
\begin{center}
%\missingfigure{multiplicative epsilon}
\includegraphics[width=\textwidth]{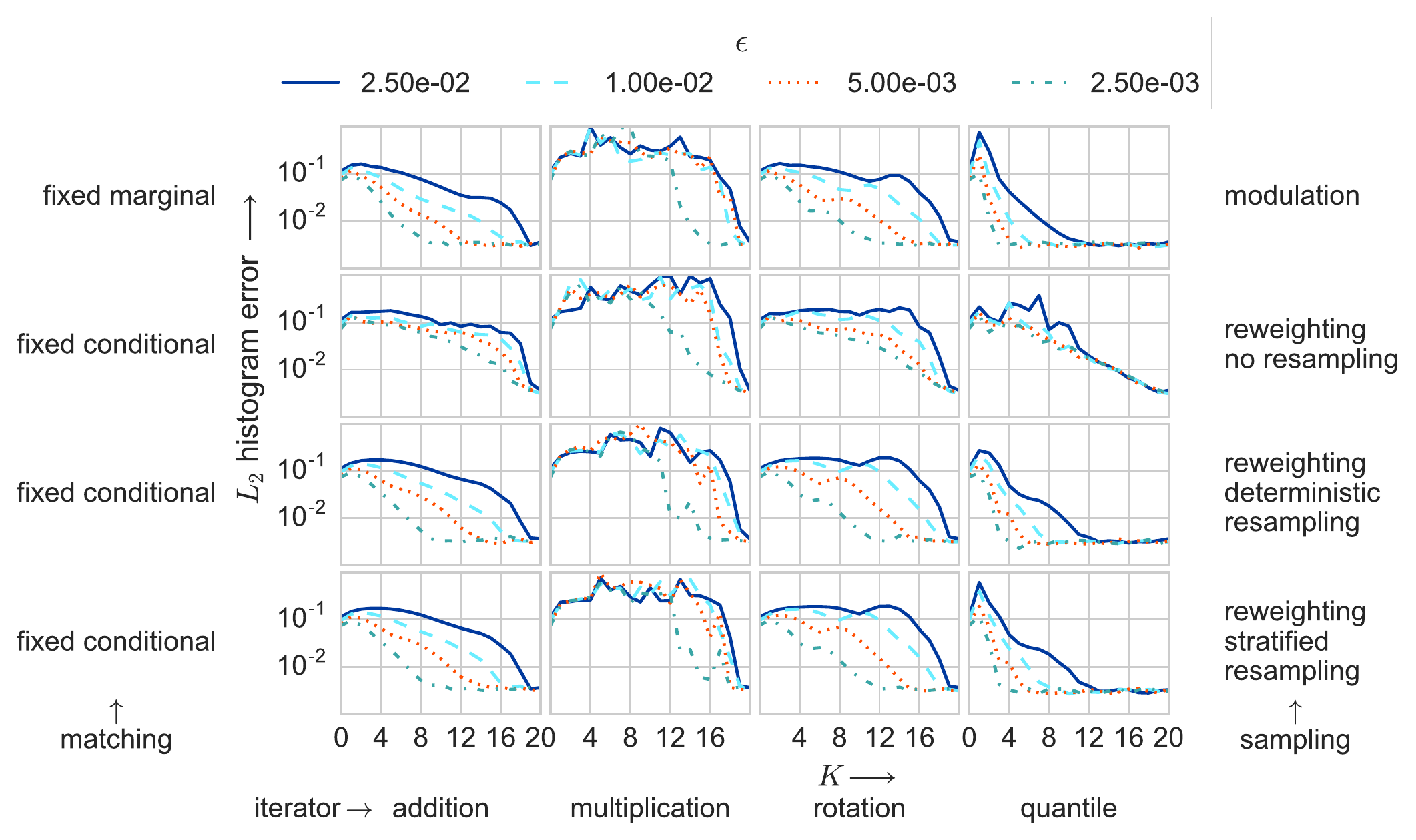}
\caption{Error in the histogram of the slow variable as a function of the number of parareal iterations in the case of a \emph{bistable} fast variable ($k_0=60$). Each graph corresponds to a different choice for the matching and sampling (rows) and for the iterator (columns). The different lines represent results for different values of $\epsilon$, indicating the time scale of the fast variable. Details of the simulations are presented in the text. The reference trajectory (used to compute errors) is a serial simulation with a tenfold increase in the number of particles.}
\label{fig:regimes}
\end{center}
\end{figure}

Figure~\ref{fig:regimes} shows the convergence behaviour for all possible combinations of matchings, samplings and iterators. 
The results are given in terms of the $L^2$ error on the histogram of the slow variable versus the number of parareal iterations $k$. Errors are computed with respect to a reference, serial solution of the microscopic dynamics with a tenfold increase in the number of particles. 
The results lead to three main observations. 
First, for moderate time scale separation, the problem is more challenging and the algorithm converges slowly for some implementations of the algorithm. 
Second, the choice of iterator is more important than the choice of matching and sampling. 
Only when using the quantile iterator does the algorithm converge steadily for the challenging case where $\epsilon = 2.5 \times 10^{-2}$. 
Finally, it is evident that reweighting without resampling is a poor choice. Regardless of the iterator, the algorithm with this sampling procedure performs worse than the other approaches. 

\subsection{Scaling with number of intervals}
%\label{sec:results:bistable}

We now focus on the most challenging case discussed in Section~\ref{sec:results:epsilon}, where $\epsilon = 2.5 \times 10^{-2}$. With the other parameters as in Table~\ref{tab:default_parameters}, this puts the current simulation in Regime 2, as described in Section~\ref{sec:mic}. This case is studied extensively by \citet{BrChSm14}, who investigate how the rate of switching between the two modes affects simulations of the fine and coarse models. They proceed to construct a more accurate coarse model, that uses knowledge on the bistable nature to adjust the effective potential appropriately. In our work, we leave the coarse model unchanged and study the effectiveness of the parareal algorithm in adjusting for the inaccurate coarse model and how this is affected by the number $N$ of parareal time intervals. In the present regime, switching does not occur frequently enough for the coarse propagator to give an accurate description of the underlying dynamics.

Since we do not want to focus on stopping criteria, the maximum number of parareal iterations is fixed at $K=20$, while the number $N$ of parareal reporting intervals is varied over the values $20$, $50$ and $100$.
Note that the first of these cases, $N=20$, leads to a guaranteed exact result after $K=20$ iterations. 
All other parameters are as in Table~\ref{tab:default_parameters}. 
The results are given in terms of $L^2$ error on the histogram of the slow variable versus the number of iterations $k$, where the error is again computed with respect to a reference, serial solution of the microscopic dynamics with a tenfold increase in the number of particles. 

\begin{figure}[hpt]
\begin{center}
\includegraphics[width=\textwidth]{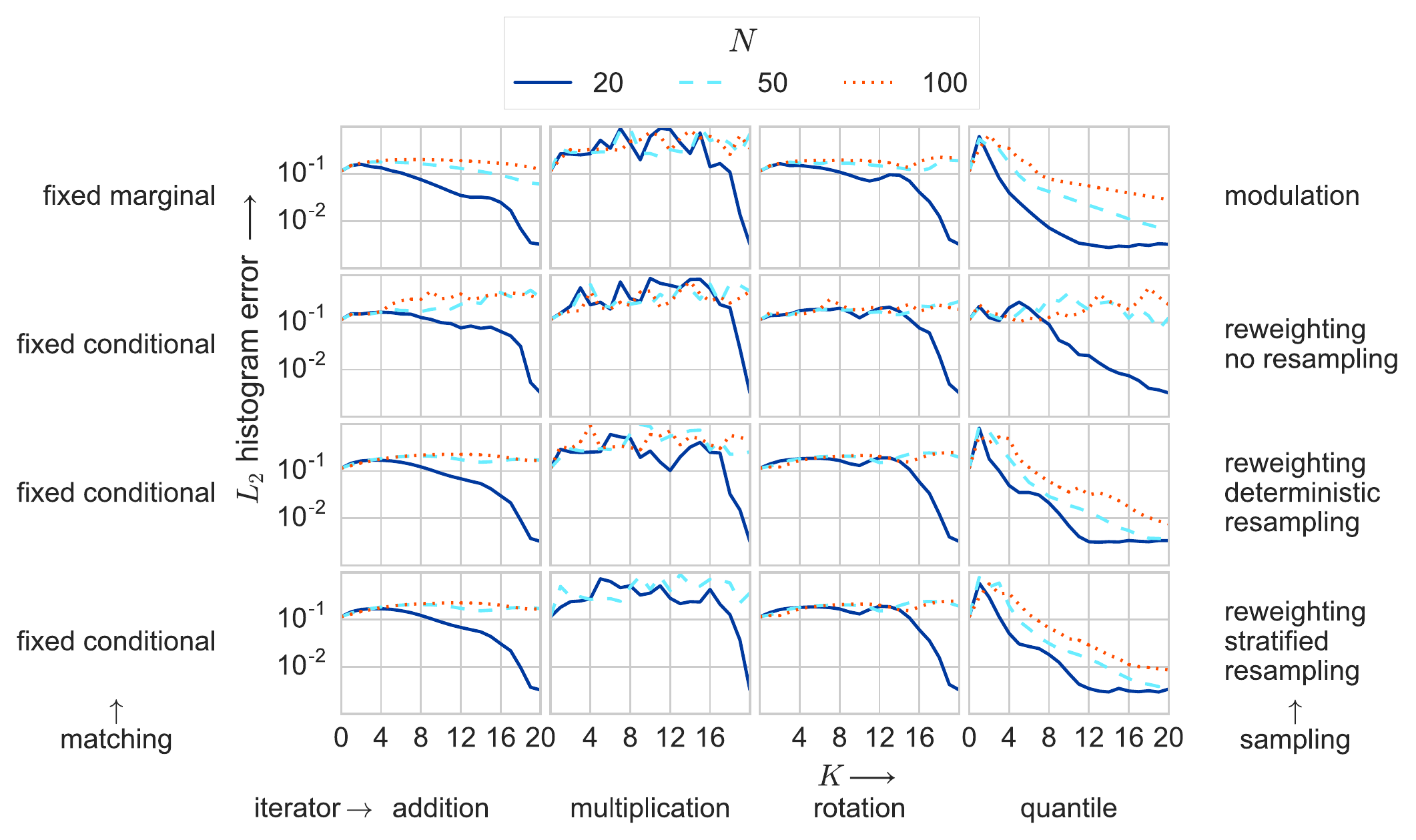}
\caption{Error in the histogram of the slow variable as a function of the number of parareal iterations in the case of a \emph{bistable} fast variable ($k_0=60$). Each graph corresponds to a different choice for the matching and sampling (rows) and for the iterator (columns). The different lines represent results for different values of $N$, the number of parareal time intervals. Details of the simulations are presented in the text.}
\label{fig:metastable}
\end{center}
\end{figure}

Figure~\ref{fig:metastable} shows the convergence behaviour for all possible combinations of matchings, samplings and iterators. The required wall clock time scales roughly as $K_{req}/N$, indicating that running the algorithm with large $N$ is likely to be advantageous, if the number of parareal iterations does not increase too much. 
The impact of the choice of iterator is more pronounced in this case, and the bimodal distribution that arises for the slow variable after some evolution appears to be challenging. Only the quantile iterator succeeds in producing results that would allow the parareal algorithm to be efficient, that is to reach convergence after a number of parareal iterations $K_{req}$ considerably smaller than the number $N$ of time intervals. 
For the other iterators, only the simulation with $N=20$ converges (up to statistical error) within the allowed number of parareal iterations. And even in this case, the number of iterations which are required is too large for the algorithm to be of practical use. As for the experiments described in Section~\ref{sec:results:epsilon}, there is also a notable difference in the performance between the different matching and sampling strategies, with the reweighting without resampling underperforming again. 
%Due to a less advective solution, however, the effect is less pronounced. 

\subsection{Linear advection case}

At this point, we have proposed a variety of matching procedures and iterators, and we have numerically observed their performance on a challenging test case. A complete understanding is however still lacking, and a detailed numerical analysis of the method would certainly be a first step in that direction.
%We have however not \emph{explained} the observations, and we do not have a full understanding of the reasons for the observed behavior. A detailed numerical analysis seems therefore to be required.
We describe now a test-case aimed at better assessing the influence of the choice of iterator. 

\begin{figure}[hpt]
\begin{center}
\includegraphics[width=\textwidth]{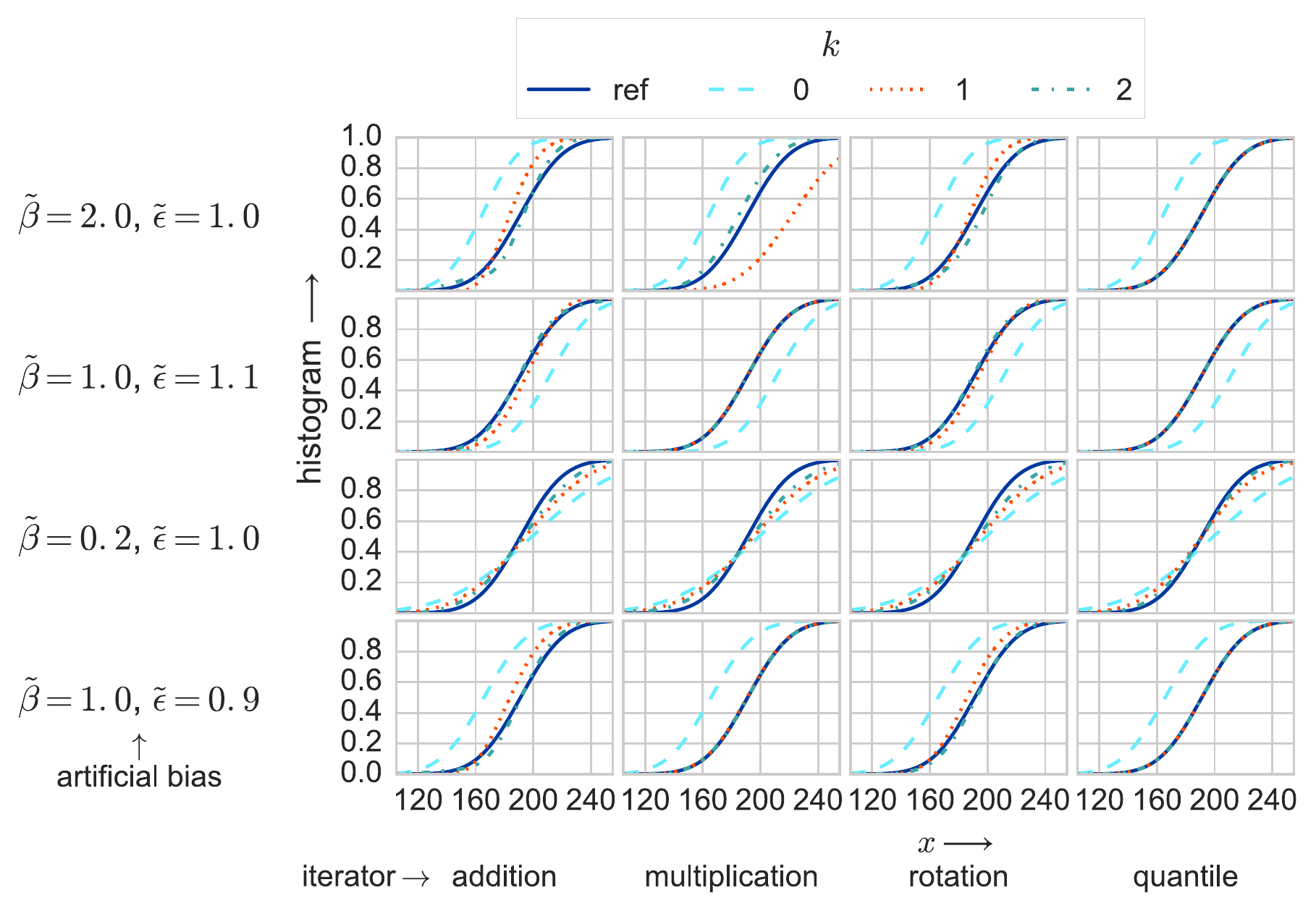}
\caption{Several iterations of the final time solution for the linear advection model. Each graph represents results for a different choice of the coarse model error (rows) and of the iterator (columns). The different lines represent different parareal iterations.}
\label{fig:dummy_k}
\end{center}
\end{figure}

We expect the superior performance of the quantile iterator to be related to the way the iterator handles differences in either the diffusive or the advective parts of the evolution between the coarse and the fine propagators. To numerically assess this intuition, we briefly discuss a different setting, where we discard the microscopic, Monte Carlo simulations and use the same discretised one-dimensional Fokker-Planck equation for both the coarse and the fine propagators. 
In this case the matching and restriction are trivial, but the iterator is not. 
This system has a potential $\overline{F}(x) = x/\epsilon$ and the diffusion is $\overline{D}(x) = 1/\beta$. We choose $\epsilon=\beta=1$. 
The coarse propagator has the same potential and diffusion, but with arbitrarily introduced biased values for the parameters: $\tilde\epsilon$ and $\tilde\beta$. Results are shown on Figure~\ref{fig:dummy_k} for different values of $\tilde\epsilon$ and $\tilde\beta$ and different choices for the iterator.

From the parareal iterations shown on Figure~\ref{fig:dummy_k}, it is obvious that the quantile iterator is capable of correcting for an error in the advection almost perfectly, as well as for an underdiffusive (top row) coarse propagator. 
Only the overdiffusive (third row) case still shows a visible error after a single iteration. 
The superior performance of the quantile iterator in this case is readily explained by the fact that the coarse propagator error is only an error in the advection and diffusion. 

\section{Conclusion}
\label{sec:conclusion}

We have presented an extension of the parareal algorithm to a multiscale setting where the fine propagator is given by a Monte Carlo simulation of a system of stochastic differential equations and the coarse  propagator models the evolution of the distribution of a scalar-valued (slow) component of the system using a discretization of the corresponding (slow) Fokker-Planck equation. 
The parareal algorithm requires a coupling between the two system descriptions. 
Finding the macroscopic probability distribution for a given microscopic ensemble of particles is trivial, but generating a microscopic ensemble consistent with a certain macroscopic distribution is more involved. 
We study several strategies (presented in Section~\ref{sec:matching}), which we conceptually divide into two parts.
The first part (matching) constructs a microscopic distribution which is consistent with the macroscopic one, and that is in some sense close to the prior distribution provided by the most recent microscopic ensemble for that time instance. 
The second part (sampling) produces an ensemble sampling this distribution. 
In practice, it is not necessary to compute the microscopic distribution explicitly. 
Provided the coupling between these two system descriptions satisfies the two properties~\eqref{eq:consistency} and~\eqref{eq:projection}, the parallel-in-time computation of the parareal algorithm leads to a significant reduction in wall clock time which is required for a given simulation. 

Besides the choice in how to couple the microscopic and macroscopic systems, there is also a choice in how to perform the parareal updates. 
The essence of the algorithm is that the discrepancy between the coarse and fine propagator results during the parallel computations is used to update the coarse model during the serial computations. 
This is usually done by simple addition and subtraction. However, since the coarse state here represents a probability distribution, the need to conserve positivity and unit norm argue for using other iterators. 
We outline four different choices in Section~\ref{sec:jumps}.

We test the algorithm for all possible combinations of matching, sampling and iterator procedures on a chemical system. 
Parameter choices in the model allow for several types of behaviour: a monostable fast variable and a bistable fast variable with a switching time that leads to a further division into three different regimes. 
For these different cases, the performance of the various implementations varies. 
However, throughout the experiments, we observe that the quantile iterator performs as well or better than the others. 
The difference between the \emph{matching} operators is less pronounced, but for the matching that leaves the conditional distribution of the fast variable unchanged it is critical that the reweighting is followed by a resampling. 
Without resampling the algorithm fails in the simplest test case, with monostable fast variable, and also underperforms in other test cases. 

\bibliographystyle{abbrvnat}
\bibliography{references}

\end{document}

%% file: macros.tex
\let\oldsqrt\sqrt
\def\sqrt{\mathpalette\DHLhksqrt}
\def\DHLhksqrt#1#2{%
\setbox0=\hbox{$#1\oldsqrt{#2\,}$}\dimen0=\ht0
\advance\dimen0-0.2\ht0
\setbox2=\hbox{\vrule height\ht0 depth -\dimen0}%
{\box0\lower0.4pt\box2}}
\renewcommand{\epsilon}{\varepsilon}

\makeatletter
\def\psls@dotteddashed{%
[ 0 \psk@dotsep CLW add \psk@dash ] 0 setdash  1 setlinecap  stroke}
\makeatother

\newcommand{\dd}{\,\mathrm{d}}

\newcommand{\bs}[1]{\boldsymbol{#1}}

\newcommand{\fd}[2]{\frac{\dd #1}{\dd #2}}

\newcommand{\expbra}[1]{\mathbb{E} \left[ #1 \right]}

\newcommand{\ensemble}[2]{\left\lbrace #1 \right\rbrace_{#2}}

\newcommand{\half}{\frac{1}{2}}

\newcommand{\IR}{\mathbb{R}}

%% file: parareal_micmac.bbl
\begin{thebibliography}{31}
\providecommand{\natexlab}[1]{#1}
\providecommand{\url}[1]{\texttt{#1}}
\expandafter\ifx\csname urlstyle\endcsname\relax
  \providecommand{\doi}[1]{doi: #1}\else
  \providecommand{\doi}{doi: \begingroup \urlstyle{rm}\Url}\fi

\bibitem[Aubanel(2011)]{Aubanel11}
E.~Aubanel.
\newblock Scheduling of tasks in the parareal algorithm.
\newblock \emph{Parallel Computing}, 37\penalty0 (3):\penalty0 172--182, 2011.

\bibitem[Bellman(1957)]{Bellman57}
R.~E. Bellman.
\newblock \emph{Dynamic programming}.
\newblock Courier Dover Publications, 1957.

\bibitem[Bruna et~al.(2014)Bruna, Chapman, and Smith]{BrChSm14}
M.~Bruna, S.~J. Chapman, and M.~J. Smith.
\newblock Model reduction for slow--fast stochastic systems with metastable
  behaviour.
\newblock \emph{The Journal of Chemical Physics}, 140\penalty0 (17):\penalty0
  174107, 2014.

\bibitem[Dai et~al.(2013)Dai, Le~Bris, Legoll, and Maday]{DaLeLeMa13}
X.~Dai, C.~Le~Bris, F.~Legoll, and Y.~Maday.
\newblock Symmetric parareal algorithms for {H}amiltonian systems.
\newblock \emph{ESAIM: Mathematical Modelling and Numerical Analysis},
  47\penalty0 (3):\penalty0 717--742, 2013.

\bibitem[Debrabant et~al.(2017)Debrabant, Samaey, and Zieliński]{DeSaZi15}
K.~Debrabant, G.~Samaey, and P.~Zieliński.
\newblock A micro-macro method for accelerating {Monte Carlo} simulation of
  stochastic differential equations.
\newblock \emph{SIAM J. Numer. Anal.}, 55\penalty0 (6):\penalty0 2745--2786,
  2017.

\bibitem[Devroye(1986)]{Devroye86}
L.~Devroye.
\newblock Sample-based non-uniform random variate generation.
\newblock In \emph{Proceedings of the 18th Winter simulation conference}, pages
  260--265. ACM, 1986.

\bibitem[Douc et~al.(2005)Douc, Capp{\'e}, and Moulines]{DoCa05}
R.~Douc, O.~Capp{\'e}, and E.~Moulines.
\newblock Comparison of resampling schemes for particle filtering.
\newblock In \emph{Proceedings of the 4th International Symposium on Image and
  Signal Processing and Analysis, 2005. ISPA 2005}, pages 64--69. IEEE, 2005.

\bibitem[E and Engquist(2003)]{EEn03}
W.~E and B.~Engquist.
\newblock The heterogeneous multiscale methods.
\newblock \emph{Communications in Mathematical Sciences}, 1\penalty0
  (1):\penalty0 87--132, 2003.

\bibitem[E et~al.(2007)E, Engquist, Li, Ren, and Vanden-Eijnden]{EEnLiReVa07}
W.~E, B.~Engquist, X.~Li, W.~Ren, and E.~Vanden-Eijnden.
\newblock Heterogeneous multiscale methods: a review.
\newblock \emph{Commun. Comput. Phys.}, 2\penalty0 (3):\penalty0 367--450,
  2007.

\bibitem[El~{M}akrini et~al.(2007)El~{M}akrini, Jourdain, and
  Leli\`evre]{MaJoLe07}
M.~El~{M}akrini, B.~Jourdain, and T.~Leli\`evre.
\newblock Diffusion {Monte Carlo} method: Numerical analysis in a simple case.
\newblock \emph{ESAIM: Mathematical Modelling and Numerical Analysis},
  41\penalty0 (2):\penalty0 189--213, 2007.

\bibitem[Gear(2001)]{Gear01}
C.~Gear.
\newblock Projective integration methods for distributions.
\newblock \emph{NEC Trans}, 130:\penalty0 1--9, 2001.

\bibitem[Givon et~al.(2004)Givon, Kupferman, and Stuart]{GiKuSt04}
D.~Givon, R.~Kupferman, and A.~Stuart.
\newblock Extracting macroscopic dynamics: model problems and algorithms.
\newblock \emph{Nonlinearity}, 17\penalty0 (6):\penalty0 R55, 2004.

\bibitem[G\"unther et~al.(2019)G\"unther, Gauger, and
  Schroder]{guenther2018pintoptim}
S.~G\"unther, N.~R. Gauger, and J.~B. Schroder.
\newblock A non-intrusive parallel-in-time approach for simultaneous
  optimization with unsteady {PDEs}.
\newblock \emph{Optimization Methods and Software}, 34\penalty0 (6):\penalty0
  1306--1321, 2019.

\bibitem[Gurrala et~al.(2015)Gurrala, Dimitrovski, Sreekanth, Simunovic, and
  Starke]{GuDiSrSiSt15}
G.~Gurrala, A.~Dimitrovski, P.~Sreekanth, S.~Simunovic, and M.~Starke.
\newblock Parareal in time for dynamic simulations of power systems.
\newblock In \emph{International Conference on Power System Transients}. IPST,
  2015.

\bibitem[Hol et~al.(2006)Hol, Schon, and Gustafsson]{HoScGu06}
J.~D. Hol, T.~B. Schon, and F.~Gustafsson.
\newblock On resampling algorithms for particle filters.
\newblock In \emph{Nonlinear Statistical Signal Processing Workshop, 2006
  IEEE}, pages 79--82. IEEE, 2006.

\bibitem[Kevrekidis and Samaey(2009)]{KeSa09}
I.~G. Kevrekidis and G.~Samaey.
\newblock Equation-free multiscale computation: Algorithms and applications.
\newblock \emph{Annual review of physical chemistry}, 60:\penalty0 321--344,
  2009.

\bibitem[Kevrekidis et~al.(2003)Kevrekidis, Gear, Hyman, Kevrekidis, Runborg,
  and Theodoropoulos]{KeGeHyKeRuTh03}
I.~G. Kevrekidis, C.~W. Gear, J.~M. Hyman, P.~G. Kevrekidis, O.~Runborg, and
  C.~Theodoropoulos.
\newblock Equation-free, coarse-grained multiscale computation: Enabling
  mocroscopic simulators to perform system-level analysis.
\newblock \emph{Communications in Mathematical Sciences}, 1\penalty0
  (4):\penalty0 715--762, 2003.

\bibitem[Kreienbuehl et~al.(2015)Kreienbuehl, Naegel, Ruprecht, Speck, Wittum,
  and Krause]{KrNaRuSpWiKr15}
A.~Kreienbuehl, A.~Naegel, D.~Ruprecht, R.~Speck, G.~Wittum, and R.~Krause.
\newblock Numerical simulation of skin transport using parareal.
\newblock \emph{Computing and visualization in science}, 17\penalty0
  (2):\penalty0 99--108, 2015.

\bibitem[Legoll and Leli\`evre(2010)]{LeLe10}
F.~Legoll and T.~Leli\`evre.
\newblock Effective dynamics using conditional expectations.
\newblock \emph{Nonlinearity}, 23\penalty0 (9):\penalty0 2131, 2010.

\bibitem[Legoll et~al.(2013)Legoll, Leli\`evre, and Samaey]{LeLeSa13}
F.~Legoll, T.~Leli\`evre, and G.~Samaey.
\newblock A micro-macro parareal algorithm: application to singularly perturbed
  ordinary differential equations.
\newblock \emph{SIAM Journal on Scientific Computing}, 35\penalty0
  (4):\penalty0 A1951--A1986, 2013.

\bibitem[Leli\`evre et~al.(2018)Leli\`evre, Samaey, and Zieliński]{LeSaZi17}
T.~Leli\`evre, G.~Samaey, and P.~Zieliński.
\newblock Analysis of a micro-macro acceleration method with minimum relative
  entropy moment matching.
\newblock \emph{arXiv preprint 1801.01740}, 2018.

\bibitem[Lions et~al.(2001)Lions, Maday, and Turinici]{LiMaTu01}
J.-L. Lions, Y.~Maday, and G.~Turinici.
\newblock A ``parareal'' in time discretization of {PDE}'s.
\newblock \emph{Comptes Rendus de l'Acad\'emie des Sciences de Paris, Series I
  Mathematics}, 332\penalty0 (7):\penalty0 661--668, 2001.

\bibitem[Loderer et~al.(2014)Loderer, Heuveline, and Lohner]{LoHeLo14}
T.~Loderer, V.~Heuveline, and R.~Lohner.
\newblock The parareal algorithm as a new approach for numerical integration of
  {ODEs} in real-time simulations in automotive industry.
\newblock \emph{PAMM}, 14\penalty0 (1):\penalty0 1027--1030, 2014.

\bibitem[Maday and Turinici(2002)]{MaTu02}
Y.~Maday and G.~Turinici.
\newblock A parallel in time approach for quantum control: the parareal
  algorithm.
\newblock In \emph{Proceedings of the 41st IEEE Conference on Decision and
  Control, Dec. 2002, Las Vegas, USA}, volume~1, pages 62--66. IEEE, 2002.

\bibitem[Pavliotis and Stuart(2008)]{PaSt08}
G.~Pavliotis and A.~Stuart.
\newblock \emph{Multiscale Methods: Averaging and Homogenization}.
\newblock Springer, 2008.

\bibitem[Paz{\'u}rikov{\'a} and Matyska(2014)]{PaMa14}
J.~Paz{\'u}rikov{\'a} and L.~Matyska.
\newblock Convergence of parareal algorithm applied on molecular dynamics
  simulations.
\newblock \emph{MEMICS 2014}, page 101, 2014.

\bibitem[Randles and Kaxiras(2014)]{RaKa14}
A.~Randles and E.~Kaxiras.
\newblock A spatio-temporal coupling method to reduce the time-to-solution of
  cardiovascular simulations.
\newblock In \emph{2014 IEEE 28th International Parallel and Distributed
  Processing Symposium}, pages 593--602. IEEE, 2014.

\bibitem[Samaddar et~al.(2010)Samaddar, Newman, and S{\'a}nchez]{SaNeSa10}
D.~Samaddar, D.~E. Newman, and R.~S{\'a}nchez.
\newblock Parallelization in time of numerical simulations of fully-developed
  plasma turbulence using the parareal algorithm.
\newblock \emph{Journal of Computational Physics}, 229\penalty0 (18):\penalty0
  6558--6573, 2010.

\bibitem[Samaey et~al.(2011)Samaey, Leli\`evre, and Legat]{SaLeLe11}
G.~Samaey, T.~Leli\`evre, and V.~Legat.
\newblock A numerical closure approach for kinetic models of polymeric fluids:
  exploring closure relations for {FENE} dumbbells.
\newblock \emph{Computers \& Fluids}, 43\penalty0 (1):\penalty0 119--133, 2011.

\bibitem[Scott(2015)]{Scott15}
D.~W. Scott.
\newblock \emph{Multivariate density estimation: theory, practice, and
  visualization}.
\newblock John Wiley \& Sons, 2015.

\bibitem[Speck et~al.(2012)Speck, Ruprecht, Krause, Emmett, Minion, Winkel, and
  Gibbon]{SpRuKrEmMiWiGi12}
R.~Speck, D.~Ruprecht, R.~Krause, M.~Emmett, M.~Minion, M.~Winkel, and
  P.~Gibbon.
\newblock A massively space-time parallel {N}-body solver.
\newblock In \emph{Proceedings of the International Conference on High
  Performance Computing, Networking, Storage and Analysis}, page~92. IEEE
  Computer Society Press, 2012.

\end{thebibliography}
